\documentclass[11pt]{article}
\usepackage[T1]{fontenc}
\usepackage{graphicx}
\usepackage{amsmath,color,amsfonts,graphics,epsfig, float}
\usepackage[colorlinks=true]{hyperref}

\newcommand{\bb}{{\beta^2 \over 2}}

\newcommand{\PP}{P^{-1}}

\begin{document}
\title{One dimensional Fokker-Planck reduced dynamics of decision making models in Computational Neuroscience}

\author{Jos\'{e} A. Carrillo$^1$, St\'ephane Cordier$^2$, Simona Mancini$^2$\\
\vspace{5pt}\small{$^{1}$ Instituci\'o Catalana de Recerca i Estudis Avan\c cats and Departament de Matem\`atiques}\\[-8pt]
\small{Universitat Aut\`onoma de Barcelona, E-08193 Bellaterra, Spain}\\[-4pt]
\small{Email: \texttt{carrillo@mat.uab.es}} \\
\small{{\it On leave from:} Department of Mathematics, Imperial
College London, London SW7 2AZ, UK.}\\
\vspace{5pt}\small{$^{2}$ F\'ed\'eration Denis Poisson (FR 2964)}\\[-8pt]
\small{Department of Mathematics (MAPMO UMR 6628)}\\[-4pt]
\small{University of Orl\'eans and CNRS, F-45067 Orl\'eans,
France}}
\date{}
\maketitle

\begin{abstract}
We study a Fokker-Planck equation modelling the firing rates of
two interacting populations of neurons. This model arises in
computational neuroscience when considering, for example, bistable
visual perception problems and is based on a stochastic
Wilson-Cowan system of differential equations. In a previous work,
the slow-fast behavior of the solution of the Fokker-Planck
equation has been highlighted. Our aim is to demonstrate that the
complexity of the model can be drastically reduced using this
slow-fast structure. In fact, we can derive a one-dimensional
Fokker-Planck equation that describes the evolution of the
solution along the so-called slow manifold. This permits to have a
direct efficient determination of the equilibrium state and its
effective potential, and thus to investigate its dependencies with
respect to various parameters of the model. It also allows to
obtain information about the time escaping behavior. The results
obtained for the reduced 1D equation are validated with those of
the original 2D equation both for equilibrium and transient
behavior.
\end{abstract}

\section{Introduction}\label{intro}

In this work, we will propose a procedure to reduce rate models
for neuron dynamics to effective one dimensional Fokker-Planck
equations. These simplified descriptions will be obtained using
the structure of the underlying stochastic dynamical system. We
will emphasize the numerical and practical performance of this
procedure coming from ideas used in the probability community
\cite{BG} for a particular model widely studied in the
computational neuroscience literature.

We will consider a simple model \cite{RBW,CRLSF,Deco} formed by
two interacting families of neurons. We assume that there is a
recurrent excitation with a higher correlation to the activity of
those neurons of the same family than those of the other while a
global inhibition on the whole ensemble is due to the background
activity. These families of neurons are modelled through the
dynamics of their firing rate equations as in the classical
Wilson-Cowan equations \cite{WC}. The synaptic connection
coefficients $w_{ij}$, representing the strength of the
interaction between population $i$ and $j$, are the elements of a
$2\times 2$ symmetric matrix $W$ given by
$$
W=\left[
\begin{array}{ccc}
w_{+}-w_{I} &\phantom{t} & w_{-}-w_{I}\\
w_{-}-w_{I} &\phantom{t} & w_{+}-w_{I}
\end{array}
\right],
$$
Here, $w_+$ is the self-excitation of each family, $w_-$ the
excitation produced on the other family, and $w_I$ the strength of
the global inhibition. The typical synaptic values considered in
these works are such that $w_-<w_I<w_+$ leading to
cross-inhibition since $w_-<w_I$ and self-excitation since
$w_I<w_+$. Let us comment that these rate models are very
simplified descriptions of interacting neuron pools, more accurate
microscopic models introducing neuron descriptions at the level of
voltage and/or conductances probability distribution can be
derived \cite{Brunel00,BH99,BW03,CTM-PNAS04,CTRAM06,CTSM-PNAS04}.

The time evolution for the firing rates $\nu_i(t)$ of the neuronal
populations $i=1,2$ as given in \cite{Deco} follows the following
Stochastic Differential Equations (SDE):
\begin{equation}\label{ODE}
\tau \frac{d \nu_i(t)}{dt} = - \nu_i(t) + \phi \left( \lambda_i +
\sum_{j=1,2} w_{ij} \nu_j(t) \right) + \beta_i  \xi_i(t), \quad
i=1,2,
\end{equation}
where $\tau = 10^{-2} s$ is a time relaxation coefficient, which
will be chosen equal to 1 in the sequel except for the numerical
results, and $\xi_i(t)$ represents a white noise of normalized
standard deviation $\beta_i>0$. In \eqref{ODE} the function
$\phi(x)$ is a sigmoid function determining the response function
of the neuron population to a mean excitation
$x(t)=\lambda_i+\sum_j w_{ij} \nu_j$:
\begin{equation*}%\label{phi}
\phi(x)=\frac{\nu_c}{1+{\rm exp} (-\alpha (x/\nu_c -1))},
\end{equation*}
where $\lambda_i$ are the external stimuli applied to each neuron
population.

We will recall in the next section that the study of the decision
making process for the previous network can be alternatively
studied by means of the evolution of a Fokker-Planck equation in
two dimensions i.e. the plane $(\nu_1,\nu_2)$. The theoretical
study of such problem (existence and uniqueness of positive
solutions) was done in \cite{CCM}. However, we will emphasize that
due to slow-fast character of the underlying dynamical system the
convergence towards the stationary state for the corresponding
two-dimensional problem is very slow leading to a kind of
metastable behavior for the transients.

Nevertheless, the 2D Fokker-Planck equation allows us to compute
real transients of the network showing this metastable behavior.
Moreover, we can derive a simplified one dimensional SDE in
Section 3 by scaling conveniently the variables. Here, we use the
spectral decomposition and the linearized  slow manifold
associated to some stable/unstable fixed point of the
deterministic dynamical system. The obtained 1D Fokker-Planck
equation leads to a simple problem to solve both theoretically for
the stationary states and numerically for the transients. In this
manner, we can reduce the dynamics on the slow manifold to the
movement of a particle in an effective 1D potential with noise. We
recover the slow-fast behavior in this 1D reduction but, due to
dimension, we can efficiently compute its numerical solution for
much larger times than the 2D. We can also directly compute an
approximation to the 2D equilibrium state by the 1D equilibrium
onto the slow manifold since in 1D, every drift derives from a
potential.

Let us mention that another approach to get an approximation of
the 2D Fokker-Planck equation by 1D Fokker-Planck reduced dynamics
has been proposed in \cite{RL}. This approach is purely local via
Taylor expansion around the bifurcation point leading to a cubic
1D effective potential. Moreover, an assumption about the scaling
of the noise term is performed to be able to close the expansion
around the bifurcation point. Our approach is valid no matter how
far we are from the bifurcation point as long as the system has
the slow-fast character and we do not assume any knowledge of the
scaling of the noise term. Moreover, we can reconstruct the full
potential not only locally at the bifurcation point. We point out
that the results of their 1D Fokker-Planck reduction are compared
to experimental data in \cite{RL} with excellent results near the
bifurcation point. A similar applied analysis of our reduced
Fokker-Planck dynamics in a system of interest in computational
nuroscience is underway \cite{CCDM}.

Section 4 is devoted to show comparisons between the 2D and the
reduced 1D Fokker-Planck equations both for the stationary states
and the transients. We demonstrate the power of this 1D reduction
in the comparison between projected marginals on each firing rate
variable and on the slow linearized manifold. Finally, section 5
is devoted to obtain information of the simulation in terms of
escaping times from a decision state and performance in the
decision taken.

\section{The two dimensional model}\label{2Dmodel}
We will illustrate all our results by numerical simulations
performed with the physiological values introduced in \cite{Deco}:
$\alpha = 4$ and $\nu_c = 20 Hz$, $\lambda_1=15 Hz$ and
$\lambda_2=\lambda_1+ \Delta \lambda$, with $\Delta \lambda=0$ for
the unbiased case and $\Delta \lambda=0.01, 0.05, 0.1$ for the
biased case. The noise parameter is chosen as $ \beta=0.1$, and
the connection coefficients are given by $w_{+} =2.35$,
$w_{I}=1.9$ and $w_{-}=1-r(w_{+}-1)/(1-r)$ with $r=0.3$.

\begin{figure}[H]
\centering{\includegraphics[width=6cm, angle=0]{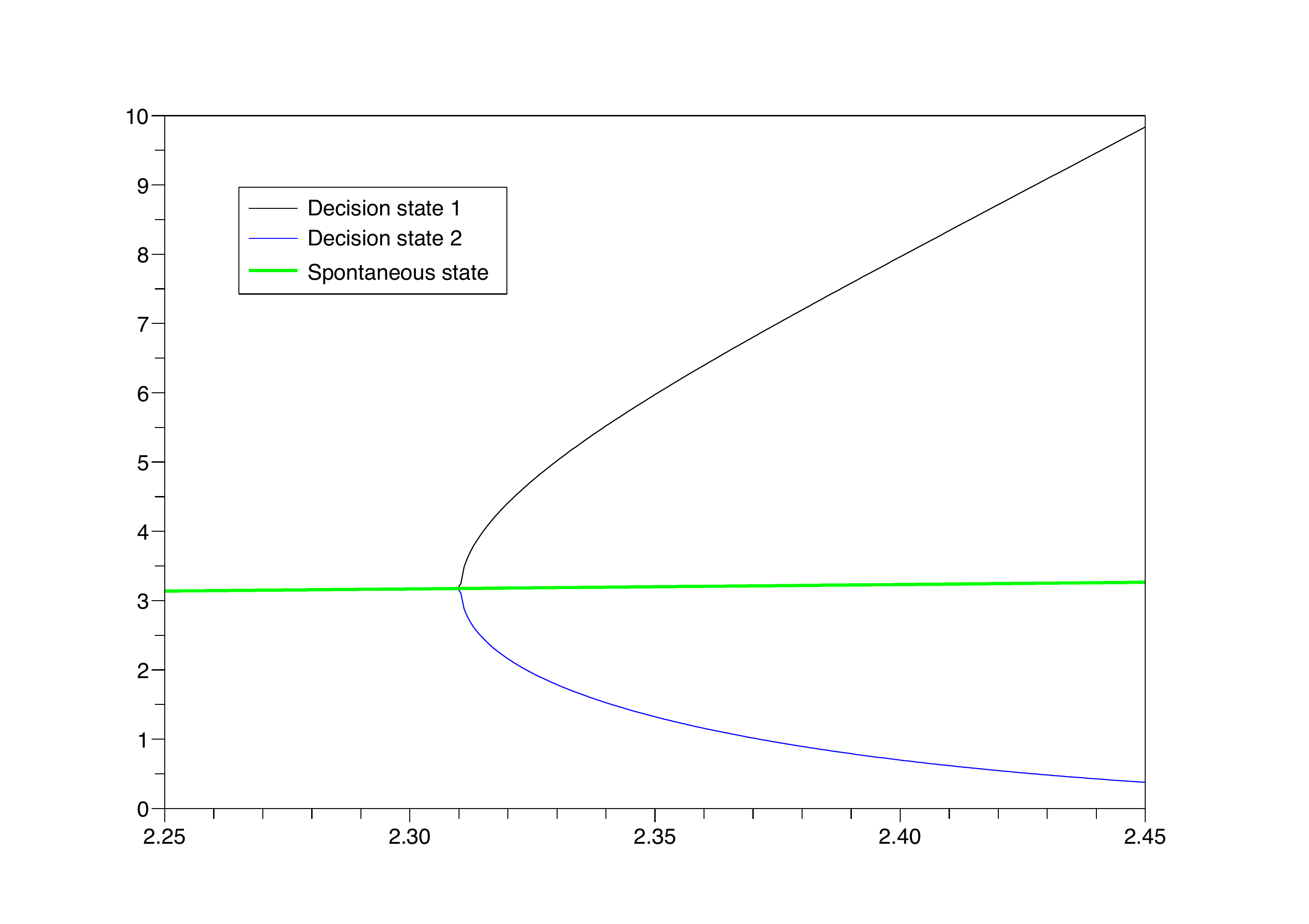}
\
\includegraphics[width=6cm, angle=0]{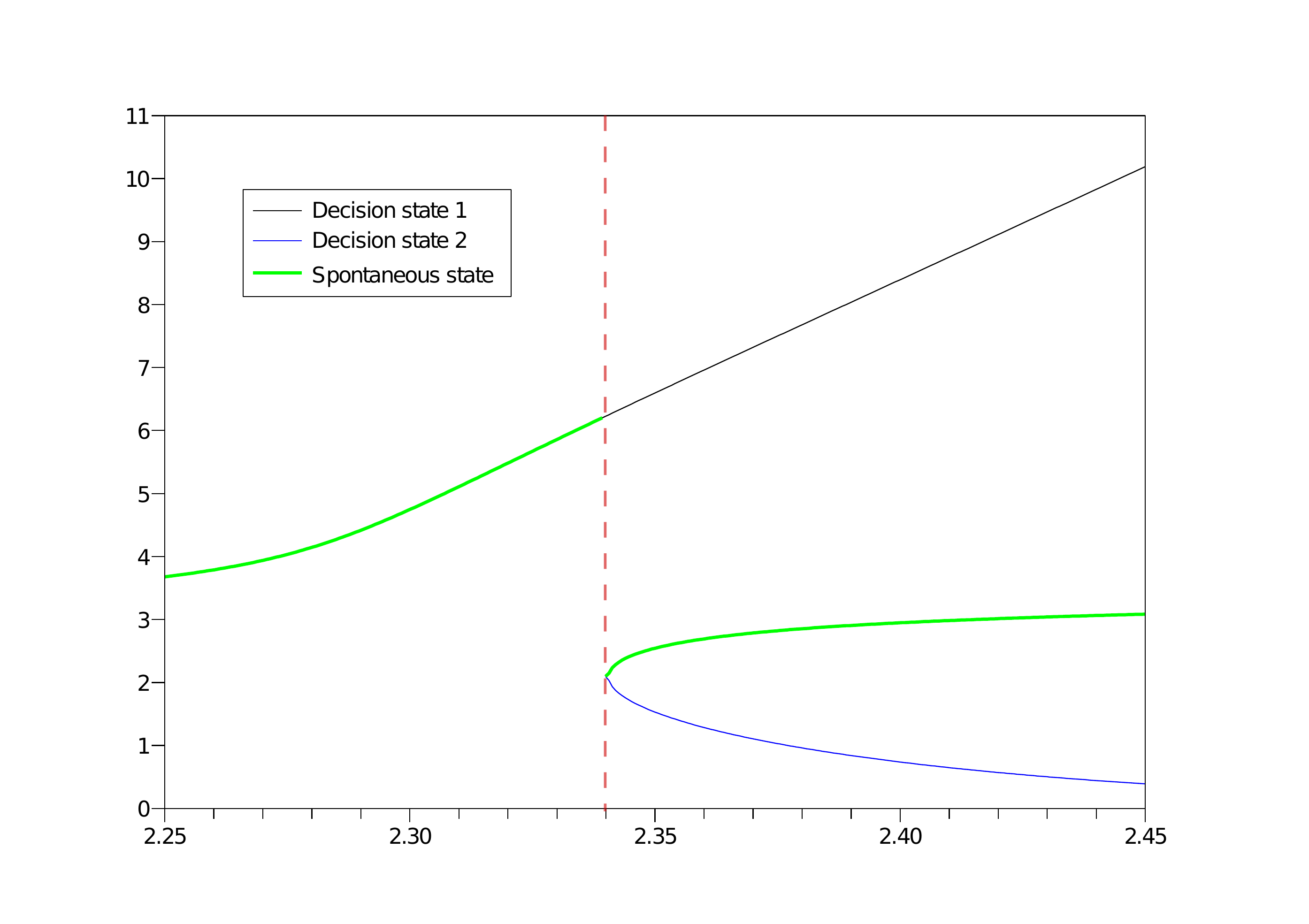}}
\caption{Bifurcation diagram: $\nu_1$-component of the equilibrium
states with respect to $w_+$. Left Figure: bifurcation diagram for
the unbiased case $\Delta \lambda =0$. Right Figure: bifurcation
diagram for the biased case $\Delta \lambda= 0.1$. }
\label{bifurcation}
\end{figure}

It is well known \cite{Deco} that the deterministic dynamical
system associated with \eqref{ODE} is characterized by a
supercritical pitchfork bifurcation in terms of the parameter
$w_+$ from a single stable asymptotic state to a two stable and
one unstable equilibrium points. We recall that the unstable point
is usually called spontaneous state while the two asymptotically
stable points are called decision states. The behavior of the
bifurcation diagram for the deterministic dynamical system
defining the equilibrium points in terms of the $w_+$ parameter
and with respect to $\Delta \lambda$ is shown in Figure
\ref{bifurcation}. Observe that in the nonsymmetric ($\Delta
\lambda\neq 0$) bifurcations, the pair of stable/unstable
equilibrium points detaches from the branch of stable points.

For example, with $w_+=2.35$ in the unbiased case, if $\Delta
\lambda =0$, the stable points are in $S_1=(1.32, 5.97)$ and its
symmetric $S_3=(5.97,1.32)$, and the unstable one is in
$S_2=(3.19, 3.19)$; whereas, in the biased case $\Delta
\lambda=0.1$ the stable points are in $S_1=(1.09, 6.59)$ and
$S_3=(5.57, 1.53)$, and the unstable one in $S_2=(3.49, 3.08)$.

\begin{figure}[H]
\centering{\includegraphics[width=6cm]{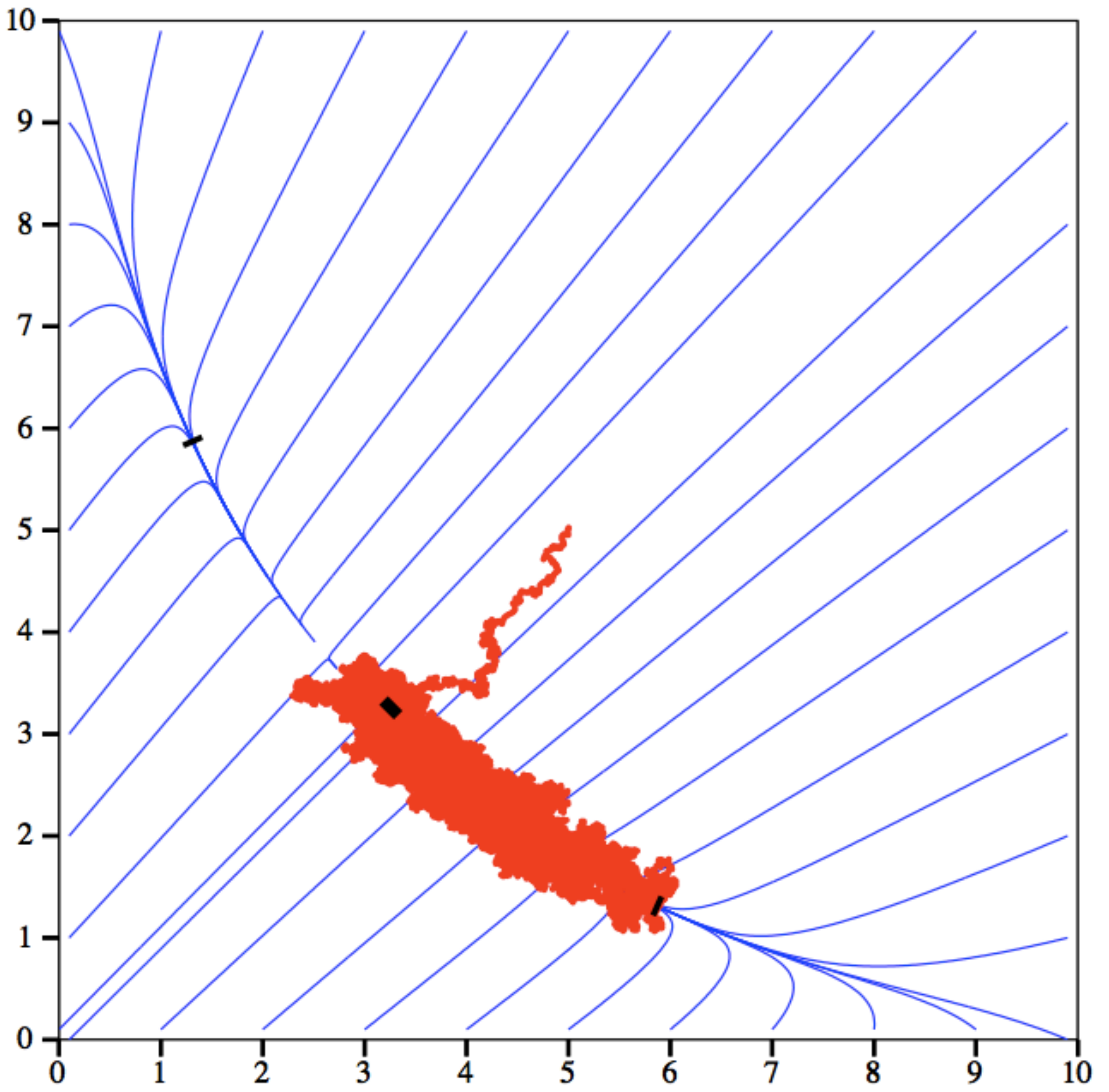}}
\caption{Dynamics of a firing rate towards stable equilibrium,
fast convergence towards the slow manifold and slow convergence
towards one of the stable equilibrium points along the slow
manifold.} \label{nils}
\end{figure}

Furthermore, it can be shown by means of direct simulations of
system \eqref{ODE}, that there is a {\it slow-fast} behaviour of
the solutions toward the equilibrium points. This behavior is
plotted in Figure \ref{nils}, where the straight lines show the
behavior of several realisations for the deterministic system
(i.e. when $\beta_i=0$), and the wiggled line represent one
realisation for the full stochastic system \eqref{ODE}. Figure
\ref{nils} highlights also the so called {\it slow manifold}: a
curves in which the three equilibrium points of the system lie and
where the dynamics are reduced to rather quickly.

Applying standard methods of Ito calculus, see \cite{Gar}, we can
prove that the probability density $p=p(t, \nu)$, with $t>0$ and
$\nu=(\nu_1, \nu_2)$, satisfies a Fokker-Planck equation (known
also as the progressive Kolmogorov equation). Hence, $p(t,\nu)$
must satisfy:
\begin{equation}\label{FPold}
\partial_t p + \nabla \cdot \left( (-\nu + \Phi(\Lambda + W \cdot \nu) )p \right) - {\bb} \Delta p =0
\end{equation}
where $\nu  \in \Omega=[0, \nu_m] \times [0,\nu_m]$,
$\Lambda=(\lambda_1,\lambda_2)$, $\Phi(x_1,x_2)=(\phi(x_1), \phi
(x_2))$, $\nabla  = (\partial_{\nu_1}, \partial_{\nu_2})$ and
$\Delta = \Delta_\nu$. We complete equation \eqref{FPold} by the
following Robin boundary conditions or no flux conditions:
\begin{equation}\label{BCold}
\left( (-\nu + \Phi)p - {\bb} \nabla p \right) \cdot n =0
\end{equation}
where $n$ is the outward unit normal to the domain.

Physically, this kind of boundary conditions means that we have no
particles incoming in the domain. This is naturally relevant for
the boundaries $[0,\nu_m]\times \{ 0\}$ and $\{ 0\} \times
[0,\nu_m]$. For the two others boundaries, $[0,\nu_m]\times \{
\nu_m\}$ and $\{ \nu_m\} \times [0,\nu_m]$, it relies on the
choice of $\nu_m$ large enough in such a way that the evolution of
our system of particles is isolated. In practice, for the choosen
parameters of our model,  $\nu_m=10$ is a
good choice.

In order to simplify notations, let us consider, from now on, the
vector field $F= (F_1,F_2)$ representing the flux in the
Fokker-Planck equation:
\begin{equation}\label{flux}
F(\nu) \stackrel{def}{=} -\nu + \Phi(\Lambda + W \cdot \nu) = \left(
\begin{array}{c}
-\nu_1 + \varphi(\lambda_1 + w_{11}\nu_1 + w_{12}\nu_2) \\
-\nu_2 + \varphi(\lambda_2 + w_{21}\nu_1 + w_{22}\nu_2)
\end{array}\right)
\end{equation}
then, equation \eqref{FPold} and boundary conditions \eqref{BCold} read:
\begin{equation}\label{FP}
\partial_t p + \nabla \cdot \left( F\, p  -  \bb \nabla p \right) =0
\end{equation}
\begin{equation}\label{BC}
\left( F\, p - \bb \nabla p \right) \cdot n =0
\end{equation}

We refer to \cite{CCM} for numerical results and a detailed
mathematical analysis of the Fokker-Planck model
\eqref{FP}-\eqref{BC}: proof of the existence, uniqueness, and
positivity of the solution, and its exponential convergence
towards the equilibrium, or stationary state. Let us just recall
that the equilibrium state cannot be analytically given because
the flux does not derive from a potential, i.e. it is not in
gradient form.

Moreover, we remark that the slow-fast structure leads to stiff
terms and thus, to small time steps and large computational time.
In fact, the slow exponential decay to equilibrium makes
impossible to wait for time evolving computations to reach the
real equilibrium. Hence, it is difficult to numerically analyze
the effect of the various parameters of the model on the
equilibrium state, and then the importance of deriving a
simplified model capable of explaining the main dynamics of the
original one is justified. Nevertheless, one could find the
equilibrium state directly by numerical methods to find
eigenfunctions of elliptic equations. The discussed slow-fast
behavior will serve us, in the sequel, to reduce the dynamics of
the system to a one dimensional Fokker-Planck equation.

\section{One dimensional reduction}\label{1Dmodel}

In this section we present the one dimensional reduction of system
\eqref{ODE}. We shall treat first the deterministic part, see
\ref{deterministic}, then the stochastic terms, section
\ref{stochastic}, and finally we describe the one dimensional
Fokker-Planck model, see \ref{FP1D}.

\subsection{Deterministic dynamical system}\label{deterministic}
The slow-fast behavior can be characterized by considering the
deterministic system of two ordinary differential equations, i.e.
(\ref{ODE}) with $\beta_i=0$. Regardless of the stability
character of the fixed point $S_2=(\nu_1^{eq},\nu_2^{eq})$, the
slow-fast behavior is characterized by a large condition number
for the Jacobian of the linearized system at the equilibrium point
$S_2$, i.e., a small ratio between the smallest and largest
eigenvalue in amplitude.

More precisely, let us write the deterministic part of the
dynamical system \eqref{ODE} as follows:
\begin{equation}\label{nudot}
\dot \nu = F(\nu),
\end{equation}
where $\nu$ is a vector and $F(\nu)=-\nu+\Phi(\Lambda + W\nu)$ is
the flux, see \eqref{flux}, as described in section \ref{2Dmodel}.
Let us denote $\nu^{eq}$ the spontaneous equilibrium point, so
that $F(\nu^{eq})=0$. By spontaneous state we mean the only
equilibrium before the bifurcation point and the unique unstable
equilibrium point after the subcritical pitchfork bifurcation.
This equilibrium point $\nu^{eq}$ is then parameterized by the
bifurcation parameter $w_+$ and it has a jump discontinuity at the
bifurcation point for nonsymmetric cases $\Delta \lambda\neq 0$.
Hence, by construction:
$$
(\nu_1^{eq},\nu_2^{eq}) = \Phi( \Lambda + W (\nu_1^{eq},\nu_2^{eq}) ).
$$
For the system (\ref{ODE}), the linearized Jacobian matrix is given by:
\begin{equation*}%\label{jac}
J_F(z_1,z_2)=\left(
\begin{array}{cc}
-1 + w_{11} \varphi'(z_1)  &  w_{12} \varphi'(z_1) \\
w_{21}  \varphi'(z_2) & -1 + w_{22} \varphi'(z_2)
\end{array}
\right) ,
\end{equation*}
where we have denoted by $z_i$ the values $z_i \stackrel{def}{=}
\lambda_i + w_{i1} \nu_1 + w_{i2} \nu_2$.

\begin{figure}[H]
\centering{\includegraphics[width=7cm]{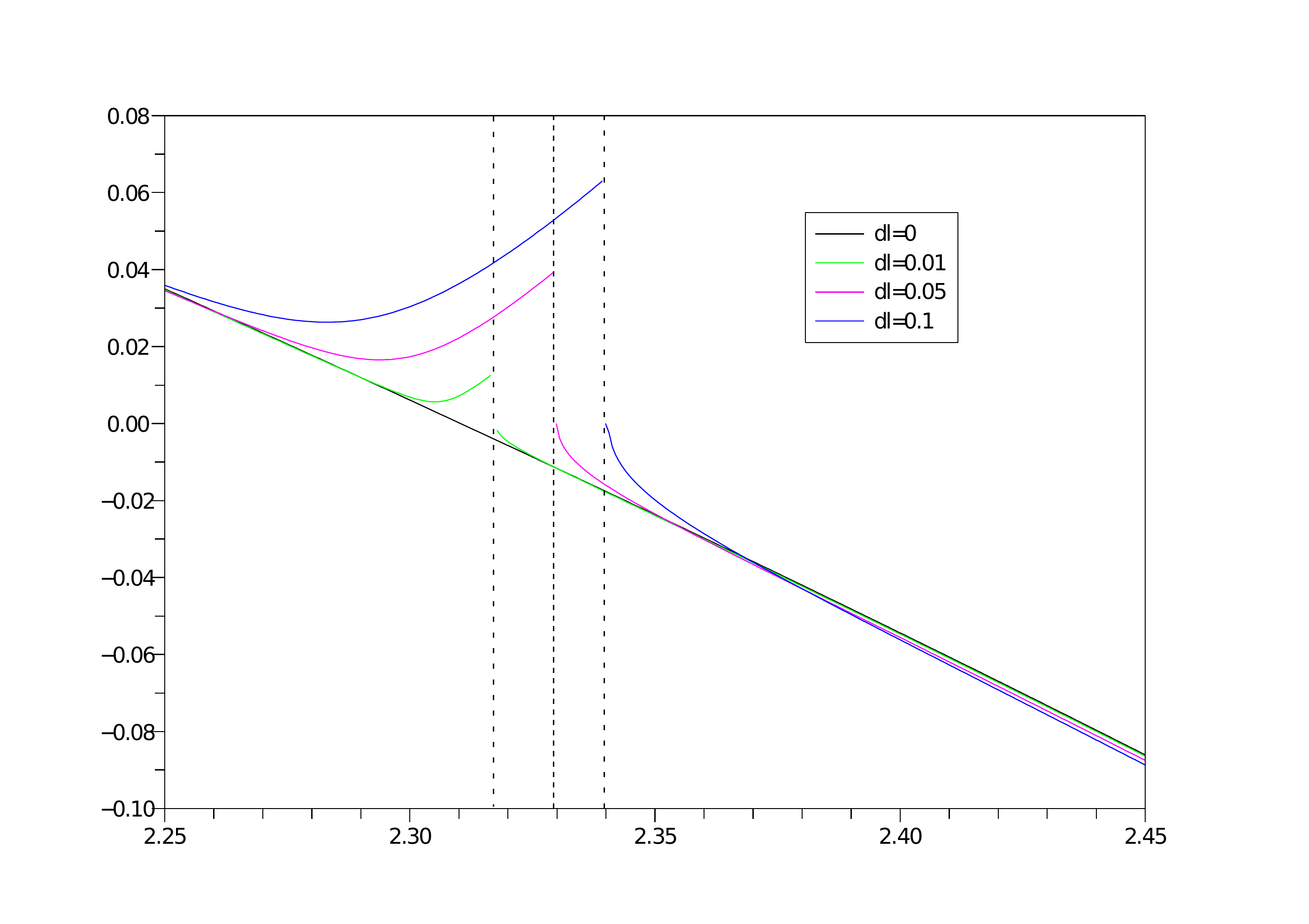}}
\caption{Eigenvalues ratio with respect to $w_+$ and for $\Delta \lambda = 0, 0.01, 0.05, 0.1$.}
\label{epsilon}
\end{figure}

We recall that $\nu^{eq}$ is an hyperbolic fixed point (saddle
point) after the bifurcation while before it is an asymptotically
stable equilibrium. Hence the Jacobian $J_F(\nu^{eq})$ has two
real eigenvalues $\mu_1$ and $\mu_2$ being both negative before
the bifurcation and of opposite signs after. The bifurcation is
characterized by the point in which the smallest in magnitude
eigenvalue becomes zero. Let us denote by $\mu_1$ the (large)
negative eigenvalue and by $\mu_2$ the (small) negative/positive
eigenvalue of $J_F(\nu^{eq})$. We remark that, the small parameter
$\varepsilon <<1$ which is responsible for the slow-fast behavior
is determined by the ratio of the amplitude of the two
eigenvalues:
\begin{equation}\label{smallpar}
\varepsilon= \left|\frac{\mu_2}{\mu_1}\right|\,.
\end{equation}
The values in terms of $w_+$ and for different values of $\Delta
\lambda$ are shown in figure \ref{epsilon}. In the range of
parameters we are interested with, $\varepsilon$ is of the order
of $10^{-2}$. Note the jump discontinuities at the bifurcation
point for $\Delta\lambda \neq 0$ since the point around which our
analysis can be performed jumps to the new created branch of the
bifurcation diagram at the bifurcation point.

\begin{figure}[H]
\centering{\includegraphics[width=5cm]{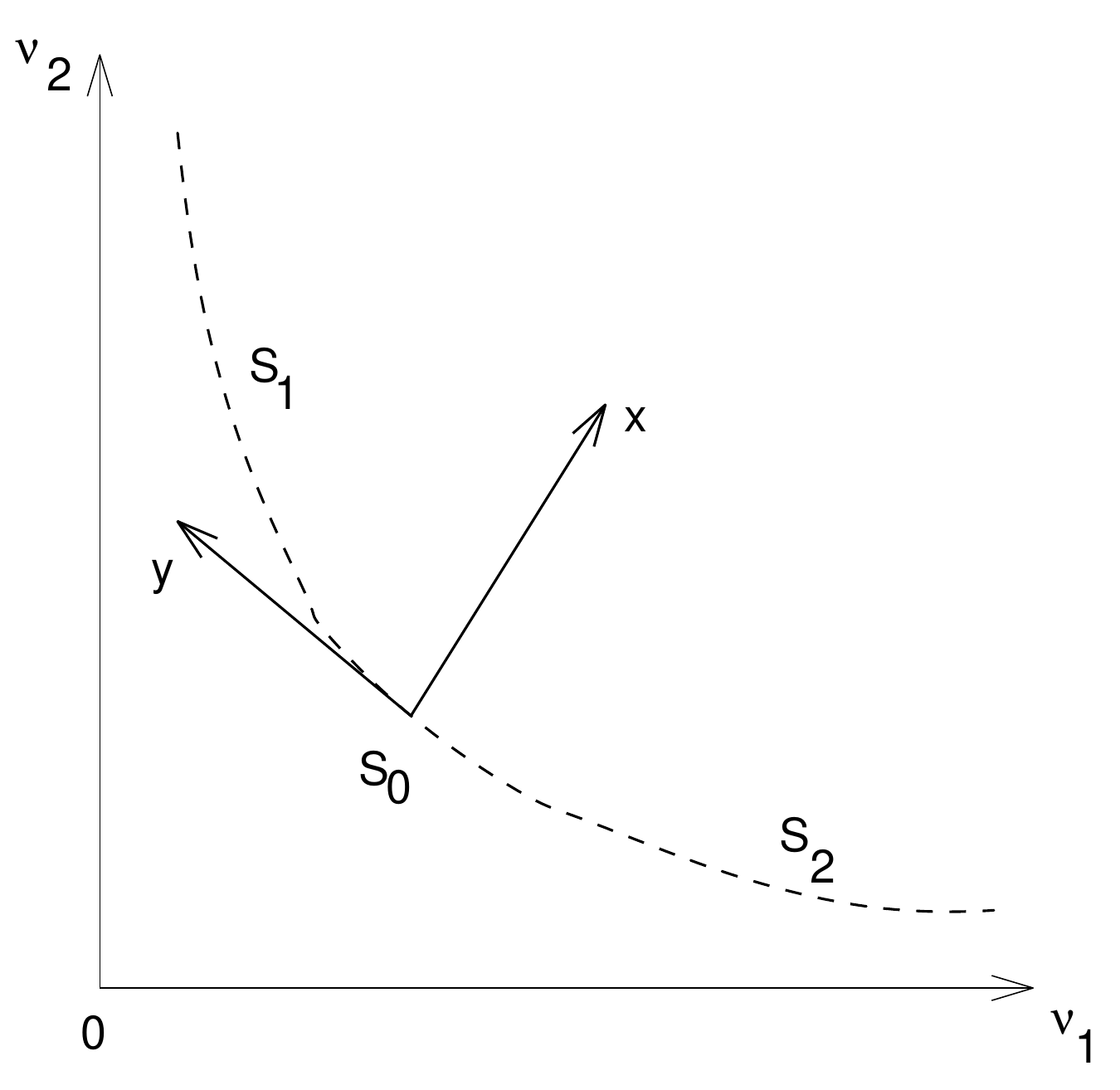}}
\caption{Change of variable from the phase space $(\nu_1, \nu_2)$ to $(x,y)$.}
\label{chang-variable}
\end{figure}

In order to reduce the system we need to introduce a new phase
space based on the linearization of the problem. We will denote by
$P$ the matrix containing the normalized eigenvectors and by $\PP$
its inverse matrix. Note that, in the unbiased case ($\Delta
\lambda =0$),  we have:
\begin{equation}\label{pppp}
P = {1 \over \sqrt{2}} \left( \begin{array}{cc}
1 & -1 \\
1 & 1
\end{array} \right) .
\end{equation}
with the associated eigenvalues $\mu_1=-1.55$ and $\mu_2=0.036$,
and the eigenvectors are orthogonal. Orthogonality of the
eigenvectors is no longer true for the nonsymmetric biased problem
$\Delta \lambda\neq 0$. Furthermore, using Hartman-Grobman theorem
\cite{Hartman,Grobman}, we know that the solutions of the
dynamical system are topologically conjugate with its
linearization in the neighbourhood of an hyperbolic fixed point,
which is valid in our case for all values of the bifurcation
parameter except at the pitchfork bifurcation. Let us write it as
follows:
\begin{equation}\label{JF}
 J_F(\nu^{eq}) = P D P^{-1},
\end{equation}
where $P$ is the matrix of eigenvectors and $D$ is the associated
diagonal matrix. We can describe the coordinates $\nu$ in the
eigenvector basis and centered on the saddle point $\nu^{eq}$ as
follows:
\begin{equation}\label{chvar}
\nu = \nu^{eq} + P X ,
\end{equation}
which gives the definition for the new variable $X=(x,y)$, see
also figure \ref{chang-variable}:
$$
X=P^{-1}(\nu-\nu^{eq}).
$$

We can conclude that system \eqref{nudot} reads in the $X$ phase space as:
\begin{equation}\label{Xdot}
\dot X = H(X)
\end{equation}
where $H(X)$ is the two dimensional vector defined by :
$$
H(X)=P^{-1} F(\nu^{eq} +PX) .
$$
We remark that by means of the chain rule, the Jacobian $J_H(X)$ is given by:
$$
J_H(X) = P^{-1} J_F(\nu^{eq} +PX) P ,
$$
and using \eqref{JF} and that $X(\nu^{eq})=0$, we obtain that
$J_H(0) = D$, which is the diagonal matrix in the change of
variables \eqref{JF}.

Let us now make explicit the system \eqref{Xdot} in terms of its
components $f(x,y)$ and $g(x,y)$:
$$\begin{cases}
\dot x = f(x,y) \\
\dot y = g(x,y)
\end{cases},
$$
where considering the definition of the flux $F$ given by
\eqref{flux}:
\begin{equation}\label{ODExy}
f(x,y)  = -x-(P^{-1}\nu^{eq})_1 + (P^{-1} \Phi (\Lambda + W
(\nu^{eq} + P X ) ))_1
\end{equation}
\begin{equation*}%\label{ODExy-g}
g(x,y)  = -y-(P^{-1}\nu^{eq})_2 + (P^{-1} \Phi (\Lambda + W (\nu^{eq} + P X ) ))_2
\end{equation*}
Now, we can choose a new time scale for the fast variable
$\tau=\varepsilon t$, with $\varepsilon$ given by
\eqref{smallpar}, in such a way that for large time $t\simeq
\varepsilon^{-1}$ then the fast variable $\tau\simeq 1$ and the
variations $\frac{dx}{d\tau}\simeq O(1)$. Then, the fast character
of the variable $x$ is clarified, see similar arguments in
\cite{BG}, and the system reads as
\begin{equation*}%\label{tildeODExy}
\begin{cases}
\varepsilon\displaystyle\frac{dx}{d\tau} = f(x,y) \\[4mm]
\displaystyle\frac{dy}{dt} = g(x,y)
\end{cases}\, .
\end{equation*}
Our model reduction assumption consists in assuming that the curve
defined by equation $f(x,y) =0$ is a good approximation when
$\epsilon \ll 1$ to the slow manifold. This manifold coincides
with the unstable manifold that joins the spontaneous point
$\nu^{eq}$ to the two other stable equilibrium points ($S_1$ and
$S_3$) after the bifurcation point while is part of the stable
manifold before the pitchfork bifurcation.

Due to the non-linearity of the function $f$, see \eqref{ODExy}
and \eqref{flux}, we cannot expect an explicit formula for
$f(x,y)=0$. Nevertheless, since $\partial_x f(0,0) \neq 0$, the
resolution in the neighborhood of the origin is given by the
implicit function theorem. Hence we can define a curve:
\begin{equation}\label{def_curve}
x=x^*(y)
\end{equation}
such that $f(x^*(y),y)=0$ in a neighbourhood of the origin. We
also note that, by construction the approximated slow manifold
$x^*(y)$, implicitly defined by (\ref{def_curve}), intersects the
exact slow manifold at all equilibrium points, i.e. where both $f$
and $g$ vanish (nullclines). Finally, we can conclude the
slow-fast ansatz, replacing the complete dynamics by the one on
the approximated slow manifold, and obtain the reduced one
dimensional differential equation:
\begin{equation*}%\label{slowODE}
\dot y = g(x^*(y),y) .
\end{equation*}

\subsection{Stochastic term}\label{stochastic}

We consider now the stochastic terms of system \eqref{ODE}. When
changing the variable form $\nu$ to $X$, also the standard
deviation of the considered Brownian motion should be modified.
Indeed the new variables $x$ and $y$ are linear combination of
$\nu_1$ and $\nu_2$. For instance, consider two stochastic
differential eqautions: $d\nu_i = \beta_i d\xi_i,$ where $\xi_i$
are two independent normalized white noises and $\beta_i$ are the
two standard deviations, and take a linear combination of $\nu_1$
and $\nu_2$ with real constant coefficients $a_1, a_2$: $x=a_1
\nu_1 + a_2 \nu_2$. Then $x$ must obey to the following stochastic
differential equation:
$$
dx=\sqrt{ (a_1 \beta_1)^2 + (a_2 \beta_2)^2} d\xi.
$$
In our case, $X = P^{-1} (\nu-\nu^{eq})$, then we have:
$$
y=(P^{-1})_{21} (\nu_1-\nu_1^{eq}) + (P^{-1})_{22} (\nu_2-\nu_2^{eq})
$$
or developing and considering $d\nu_i =\beta_i d\xi_i$,
$$
dy = (P^{-1})_{21}\beta_1 d\xi_1 + (P^{-1})_{22}\beta_2 d\xi_2
$$
Since in our model $\beta_1=\beta_2=\beta$, and considering the
above discussion, we can write for a white noise $d\xi$:
$$
dy = \beta \sqrt{ \left( (P^{-1})_{21}\right)^2 + \left( (P^{-1})_{22}\right)^2} d\xi .
$$
Finally, we conclude that the reduced one dimensional model reads:
\begin{equation}\label{ODEy}
\dot y = g(x^*(y),y) + \beta_y d\xi
\end{equation}
with $\beta_y=\beta \sqrt{ \left( (P^{-1})_{21}\right)^2 + \left(
(P^{-1})_{22}\right)^2}$. We note that in the unbiased case
$\beta_y=\beta$, since $P$ is given by \eqref{pppp}.

\subsection{One dimensional Fokker-Planck model}\label{FP1D}

We can now consider the Fokker-Planck (or progressive Kolmogorov)
equation associated to the one dimensional stochatic differential
equation \eqref{ODEy}. This gives the reduced dynamics on the
approximated slow manifold $x=x^*(y)$. Let us remark that this
reduced SDE is obtained except at the bifurcation point and
therefore valid whenever the slow-fast decomposition is verified
or in other words whenever $\epsilon$ is small.

Consider the probability distribution function $q(t,y)$, for
$t\geq 0$ and $y\in [-y_m,+y_m]$, then it must obey to the
following one dimensional Fokker-Planck equation:
\begin{equation}\label{FPy}
\partial_t q + \partial_y \left(  g(x^*(y),y) q -  {\beta_y^2 \over 2} \partial_y q \right)  =0
\end{equation}
with no-flux boundary conditions on $y=\pm y_m$:
$$
g(x^*(y),y) q -  {\beta_y^2 \over 2} \partial_y q =0\, .
$$

Since equation \eqref{FPy} is one dimensional, it is always
possible to find the effective potential $G(y)$ being the
derivative of the flux term $g(x^*(y),y)$. In other words, we can
always define the potential function:
\begin{equation*}%\label{defG}
G(y) = \int_0^y g(x^*(z),z) dz .
\end{equation*}

\begin{figure}[H]
\centering{
\includegraphics[width=6cm, angle=0]{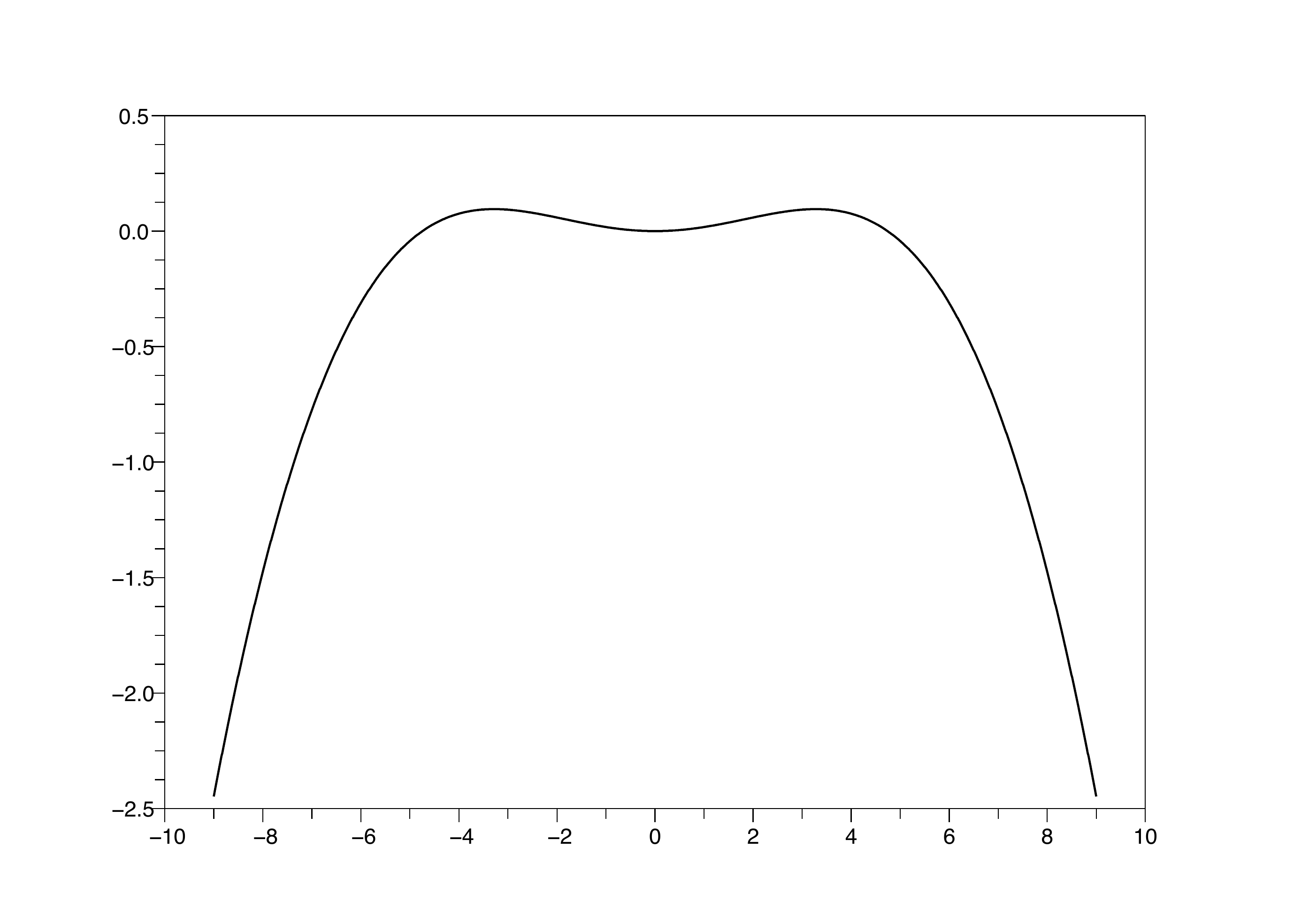}\
\includegraphics[width=6cm, angle=-0]{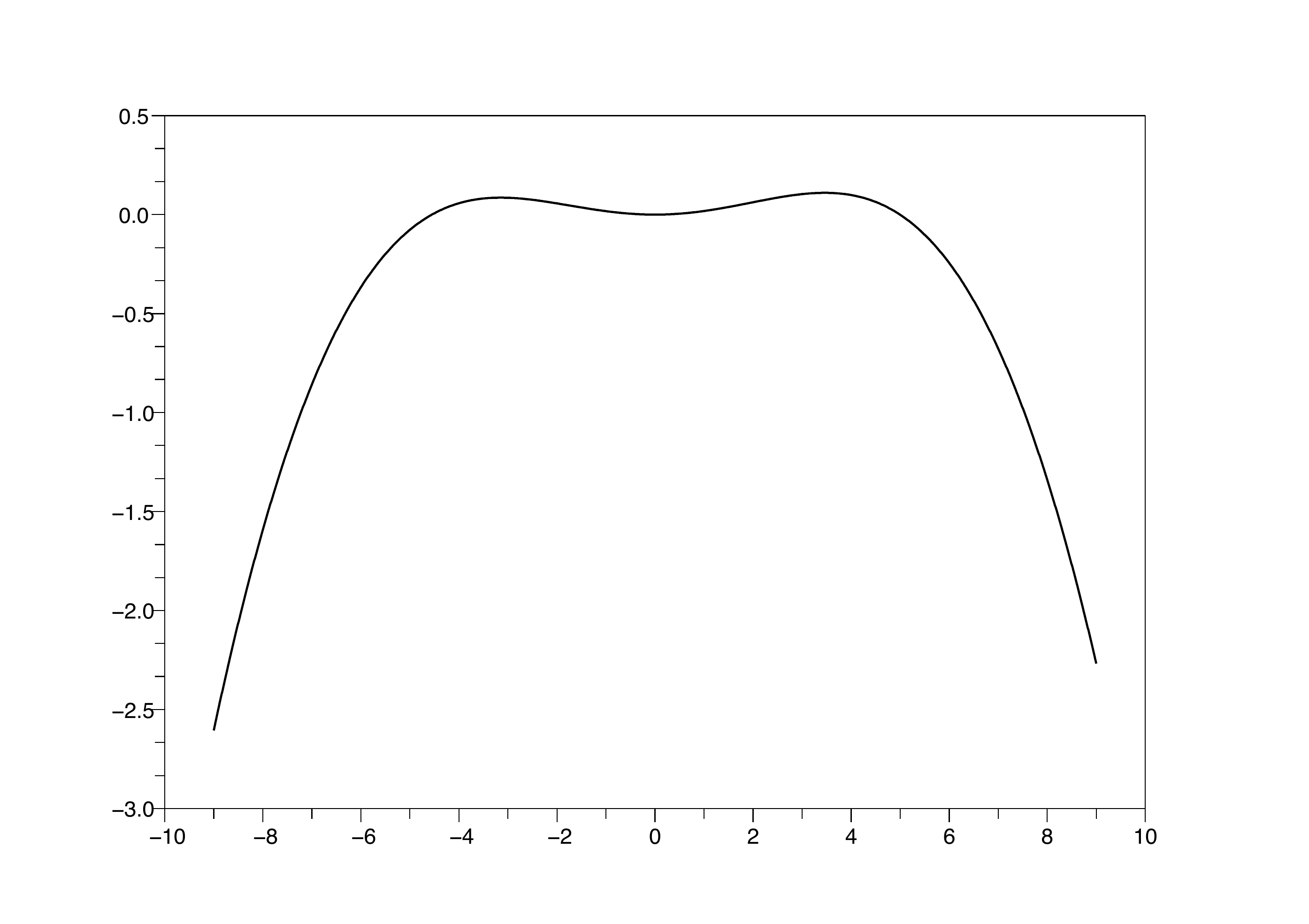} \\
\includegraphics[width=6cm, angle=-0]{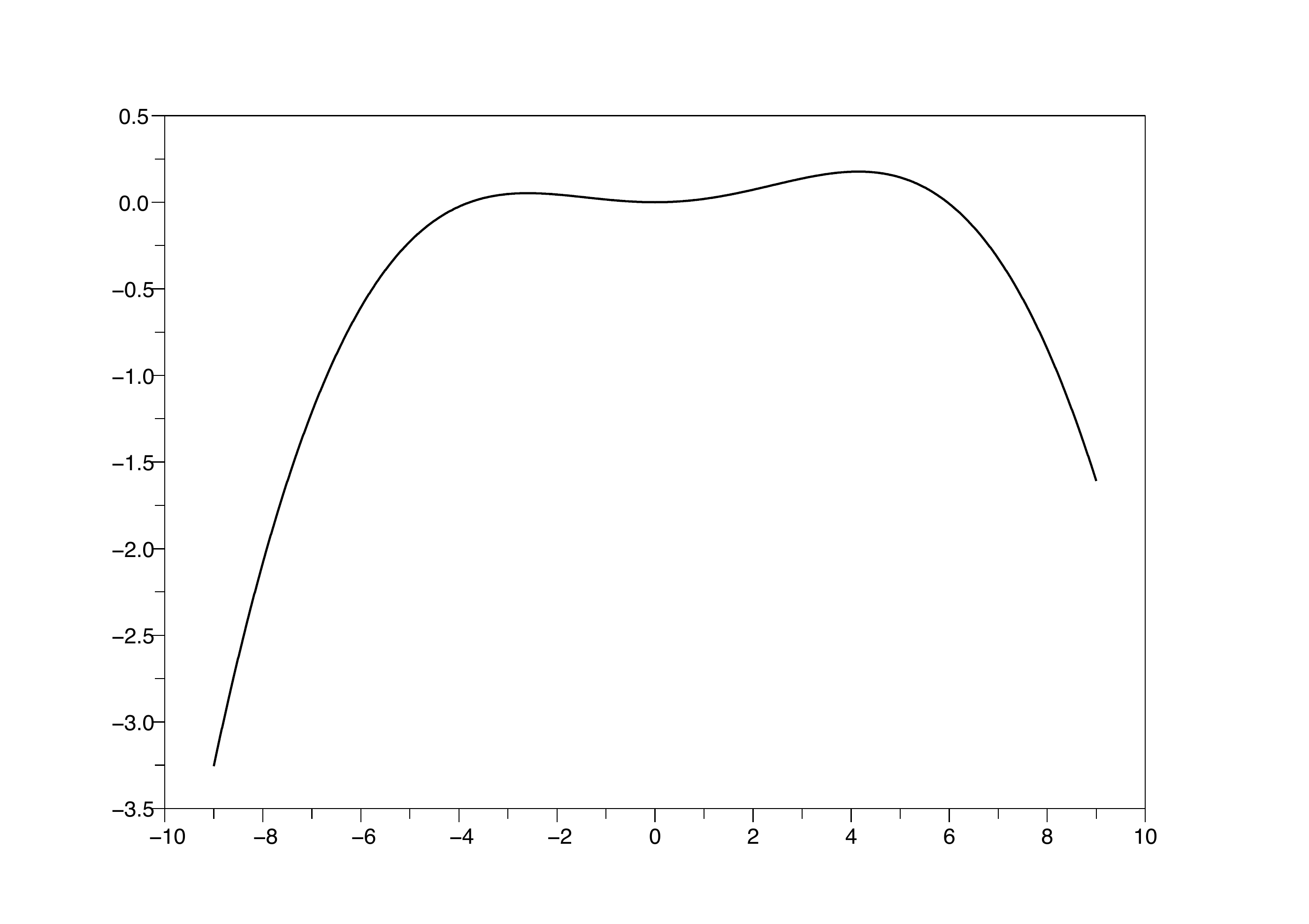}\
\includegraphics[width=6cm, angle=-0]{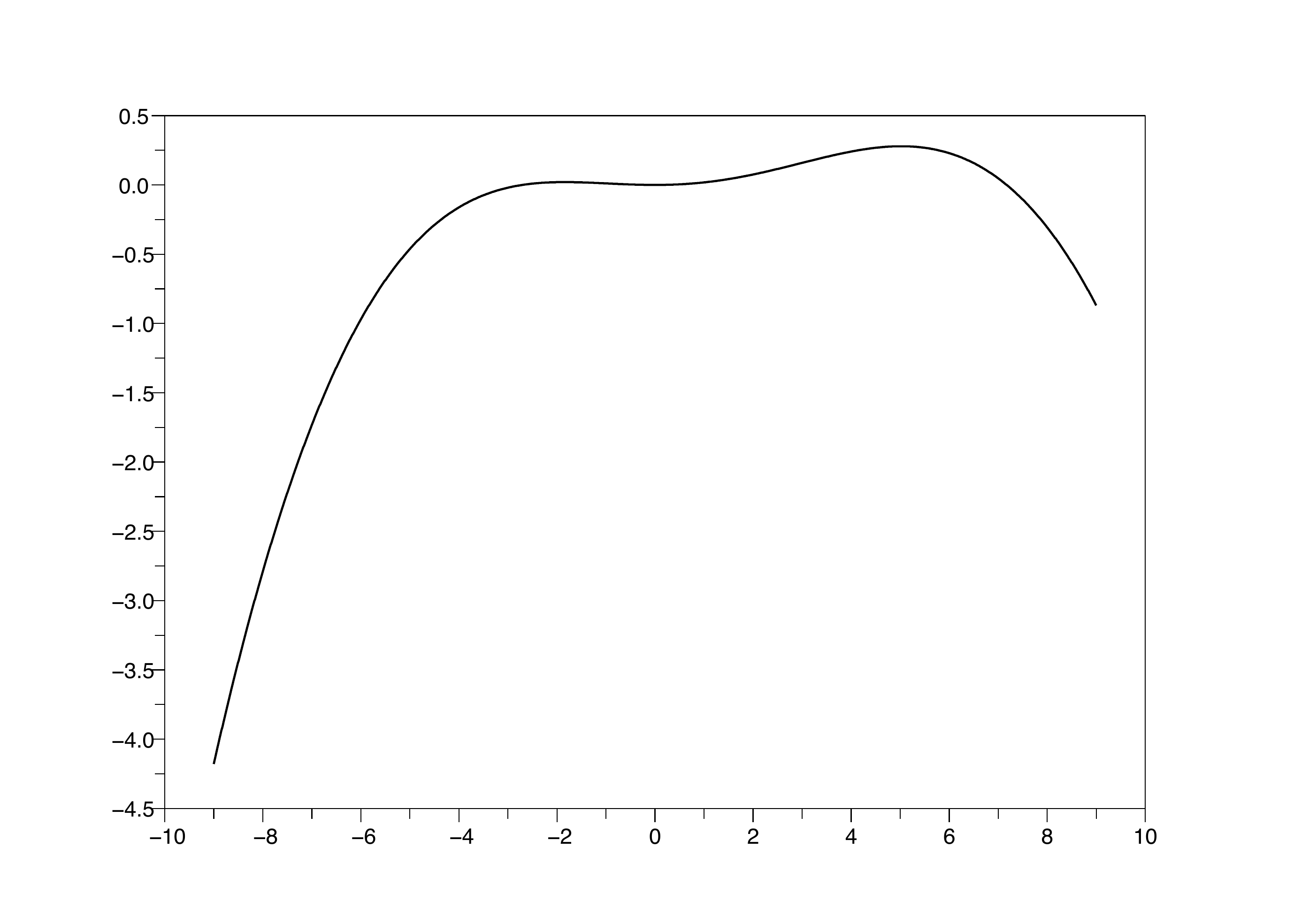}
} \caption{Comparison of the potential $-G$ for various values of
$\Delta \lambda= 0 (top-left), 0.01 (top-right), 0.05
(bottom-left), 0.1(bottom-right)$.} \label{comparison-gg}
\end{figure}

Moreover, we can explicitly obtain the stationary solutions of
\eqref{FPy}, i.e. the solutions $q_s(y)$ independent on time $t$,
as follows:
\begin{equation}\label{qstat}
q_s(y)=\frac1{Z}\exp(-2G(y)/\beta_y^2)\,,
\end{equation}
with $Z$ a suitable normalization constant. As explained also in
\cite{CCM}, these stationary solutions are the asymptotic
equilibrium states for the solution of the Fokker-Planck equation.
In other words, letting time to go to infinity, the solution
$q(t,y)$ to \eqref{FPy} must converge to $q_s(y)$. We have shown
in \cite{CCM} that the decay to equilibrium for the two
dimensional problem was exponential. Nevertheless, this decay is
so slow due to the small positive eigenvalue associated to the
spontaneous state that the simulation shows metastable behavior
for large times. Hence it is relevant to have a simple
approximated computation of their asymptotic behavior without need
to solve the whole 2D Fokker-Planck equation which is provided by
this effective 1D potential.

\section{1D model vs. 2D model}

\begin{figure}[H]
\centering{
\includegraphics[width=6cm, angle=0]{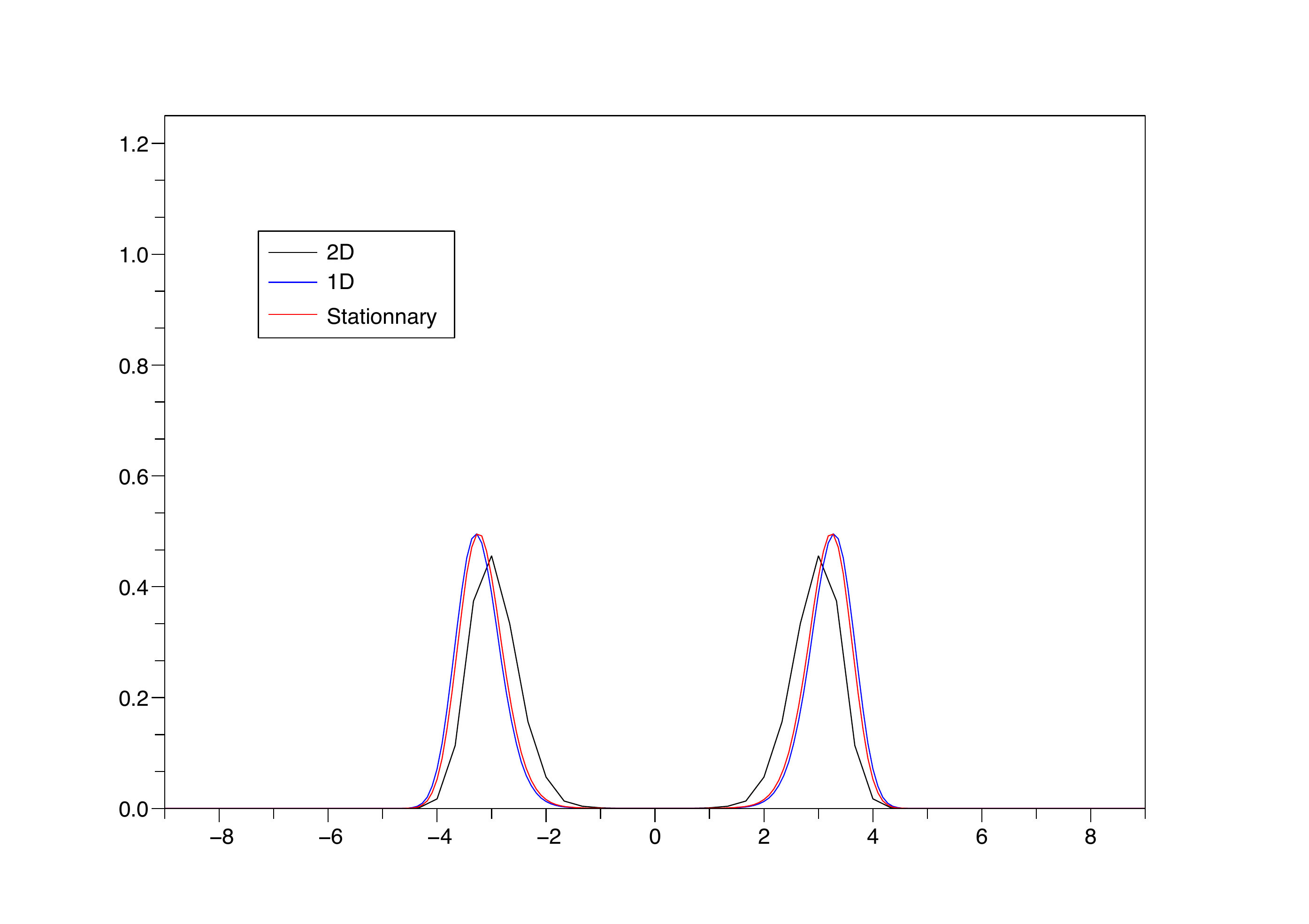}\
\includegraphics[width=6cm, angle=-0]{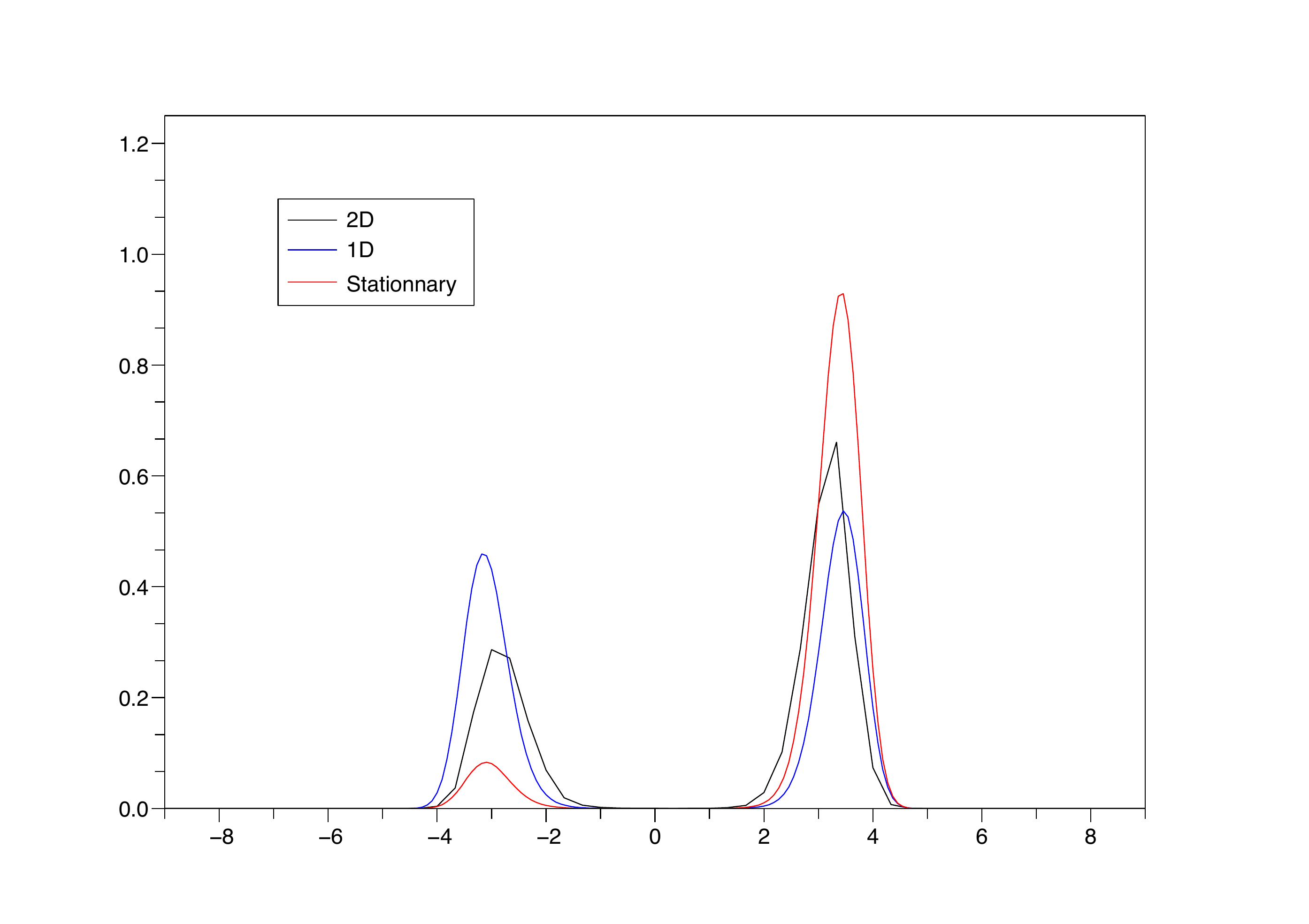} \\
\includegraphics[width=6cm, angle=-0]{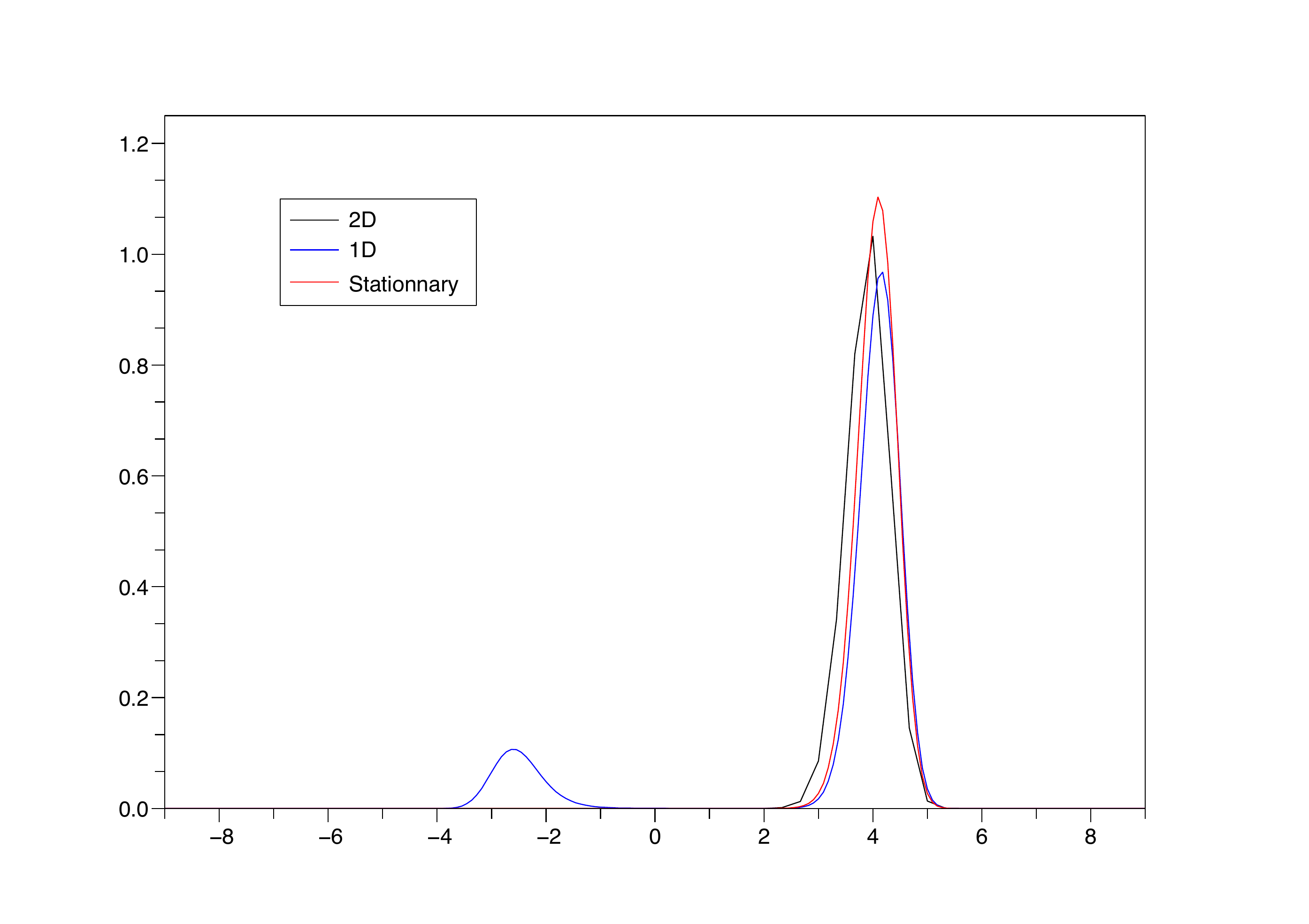}\
\includegraphics[width=6cm, angle=-0]{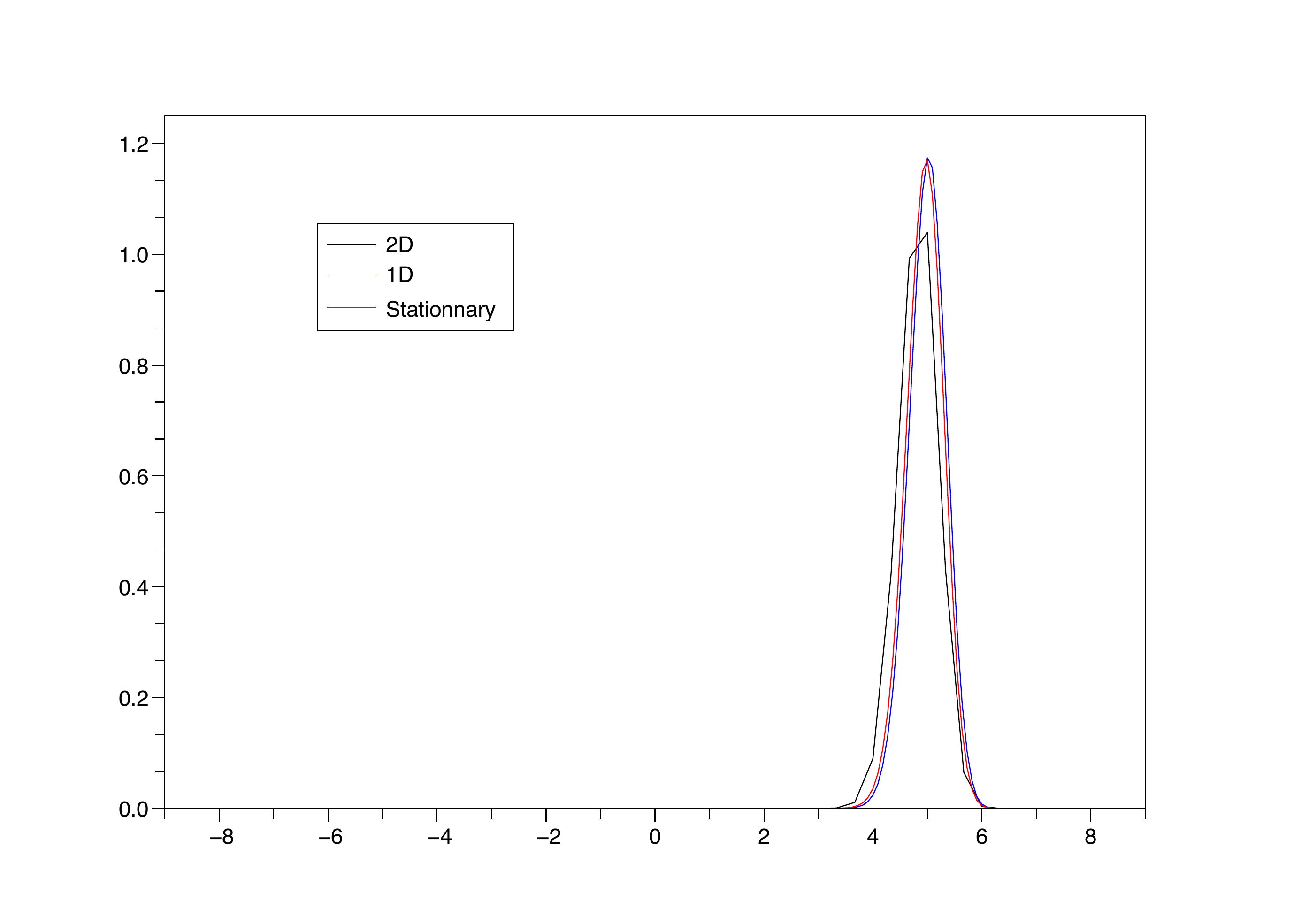}
} \caption{Comparison of the marginals in the new variable $y$,
for differents values of the biasing paramter : $\Delta \lambda= 0
(top-left), 0.01 (top-right), 0.05 (bottom-left),
0.1(bottom-right)$. Blue line: the marginal computed by means of
the 1D problem. Black line: marginal computed from the 2D problem.
Red line: the stationary marginal for the 1D. Final time is 400
seconds.} \label{comparison}
\end{figure}

In this section, we numerically compare the solutions obtained for
the one dimensional reduced Fokker-Planck equation \eqref{FPy} to
the one of the original two dimensional model \eqref{FP}.
Concerning the numerical scheme for the two dimensional problem,
we refer the reader to the detailed description in \cite{CCM}. In
particular, we are interested in the solutions at equilibrium. As
announced in section \ref{FP1D}, we have an explicit formula for
the solution at equilibrium in 1D \eqref{qstat} by computing the
primitive $G(y)$. On the contrary, in the 2D setting, we cannot
have such formulae and the computational time to approach
equilibrium is very large.

\begin{figure}[H]
\centering{
\includegraphics[width=6cm, angle=0]{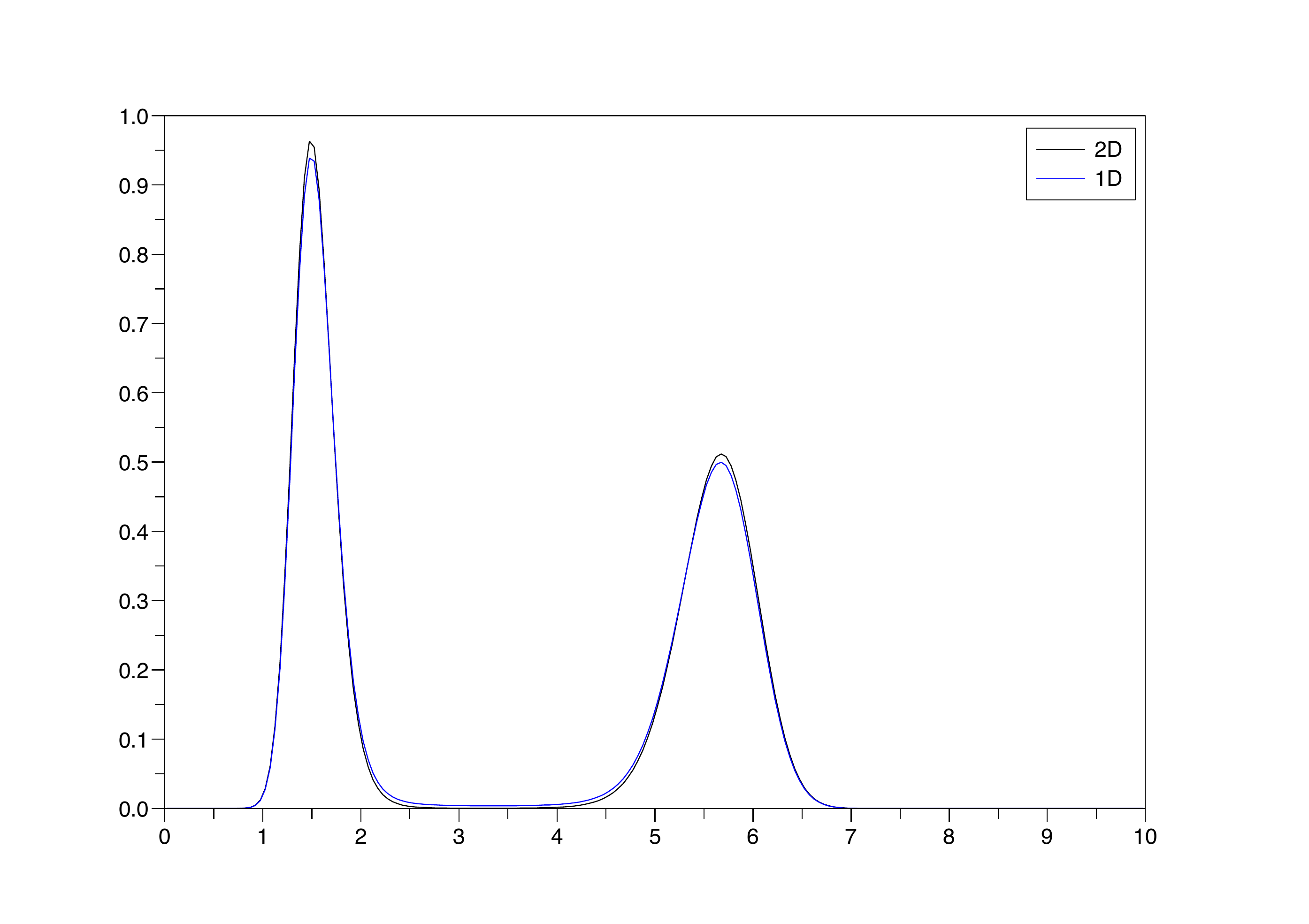}\
\includegraphics[width=6cm, angle=-0]{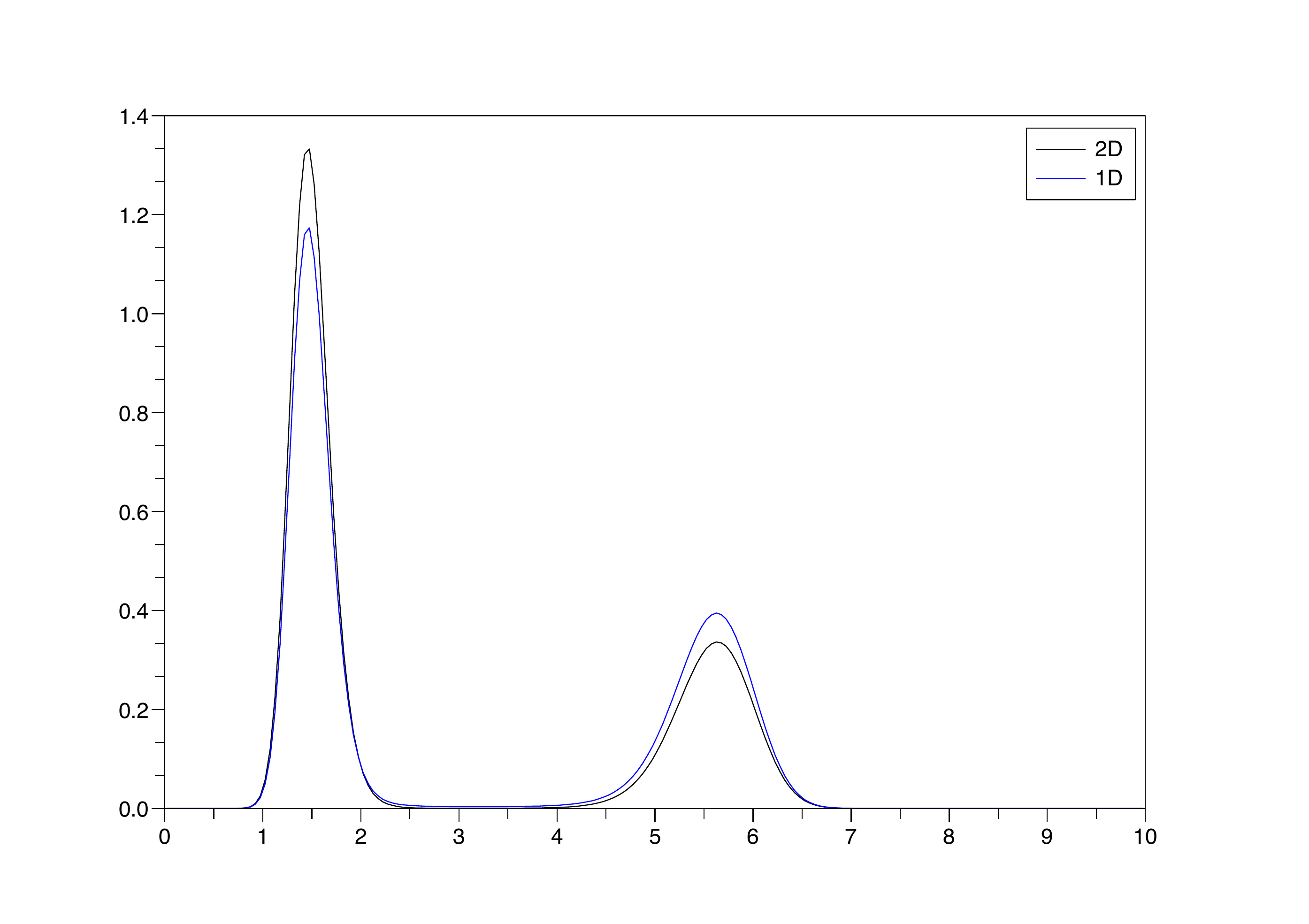} \\
\includegraphics[width=6cm, angle=-0]{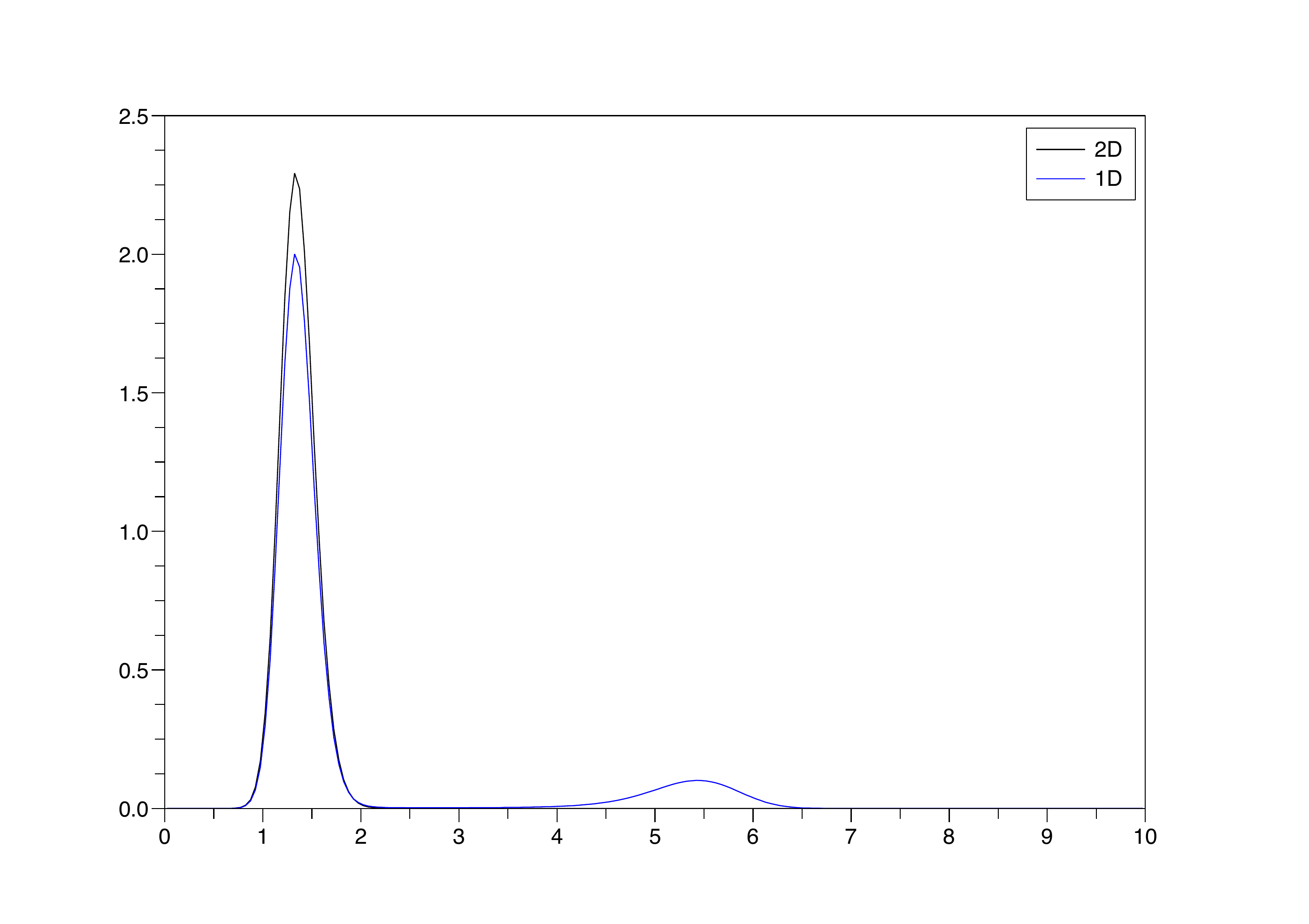}\
\includegraphics[width=6cm, angle=-0]{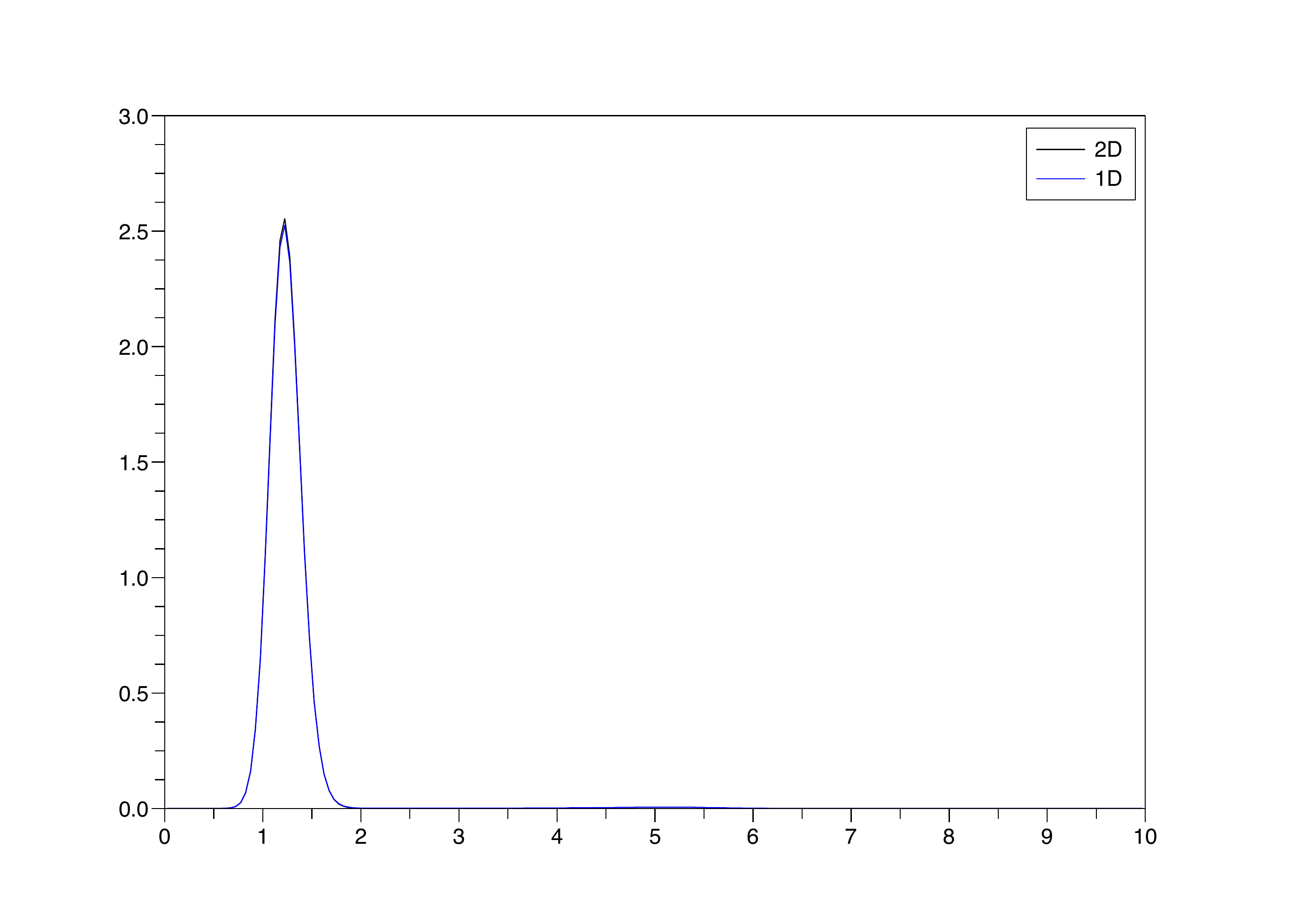}
} \caption{Comparison of the marginals in the variable $\nu_1$ for
different values of the biasing parameter: $\Delta \lambda= 0
(top-left), 0.01 (top-right), 0.05 (bottom-left),
0.1(bottom-right)$. Blue line: the marginal computed for  the 1D
problem at equilibrium. Black line: marginal computed from the 2D
problem for time $T=400 s$.  } \label{comparison2}
\end{figure}

In Figure \ref{comparison}, we plot the solution at equilibrium of
the 1D problem (the blue line) and compare it with the projection
of the two dimensional solution on the new variable $y$ (the black
line). We remark that the black line is not too smooth since we
are projecting a 2D distribution on a uniform quadrangular mesh
onto an inclined straight line. We can see a good matching in the
unbiased case ($\Delta \lambda=0$). In the biased cases, the
results are different: for $\Delta \lambda =0.01$, one clearly
sees that even if we have computed until the final time of $400$
seconds, both the 2D and the 1D solutions have not reached
equilibrium and the 2D results are closer to equilibrium; while
for $\Delta \lambda = 0.05,$ or $0.1$, the difference is smaller
since the drift is strong enough to push all particles toward only
one of the equilibrium points and there is only one population
bump at least for the 2D results. The 2D results are closer to
equilibrium at $\Delta \lambda = 0.05$ while at $\Delta \lambda =
0.1$ the 1D are closer.

On the other way round, we can also compare the marginals obtained
from the two dimensional problem with the projections of the
solution for the one dimensional problem on the $\nu_1$ and/or
$\nu_2$ axes. In figure \ref{comparison2} we show the comparison
for various $\Delta \lambda$. Note that $\Delta \lambda=0.01$ is
the most interesting case as discussed in the previous figure. In
fact, for larger $\Delta \lambda$, at equilibrium, the particles
are almost all concentrated around one of the two stable points.
Thus, no bump is visible around the second one (even in the one
dimensional reduced solution), and for the unbiased case the
matching is excellent. We warn the reader in order to compare
Figures \ref{comparison} and \ref{comparison2} that increasing
values of $y$ correspond to decreasing values of $\nu_1$.

The results demonstrate that the 1D reduction is worth to obtain
information about the behavior at equilibrium. In the next
section, we shall investigate the time dependent solution $q(y,t)$
of equation \eqref{FPy}.

\subsection{Time dependent solution}\label{time-sol}

\begin{figure}[H]
\centering{
\includegraphics[width=6cm, angle=0]{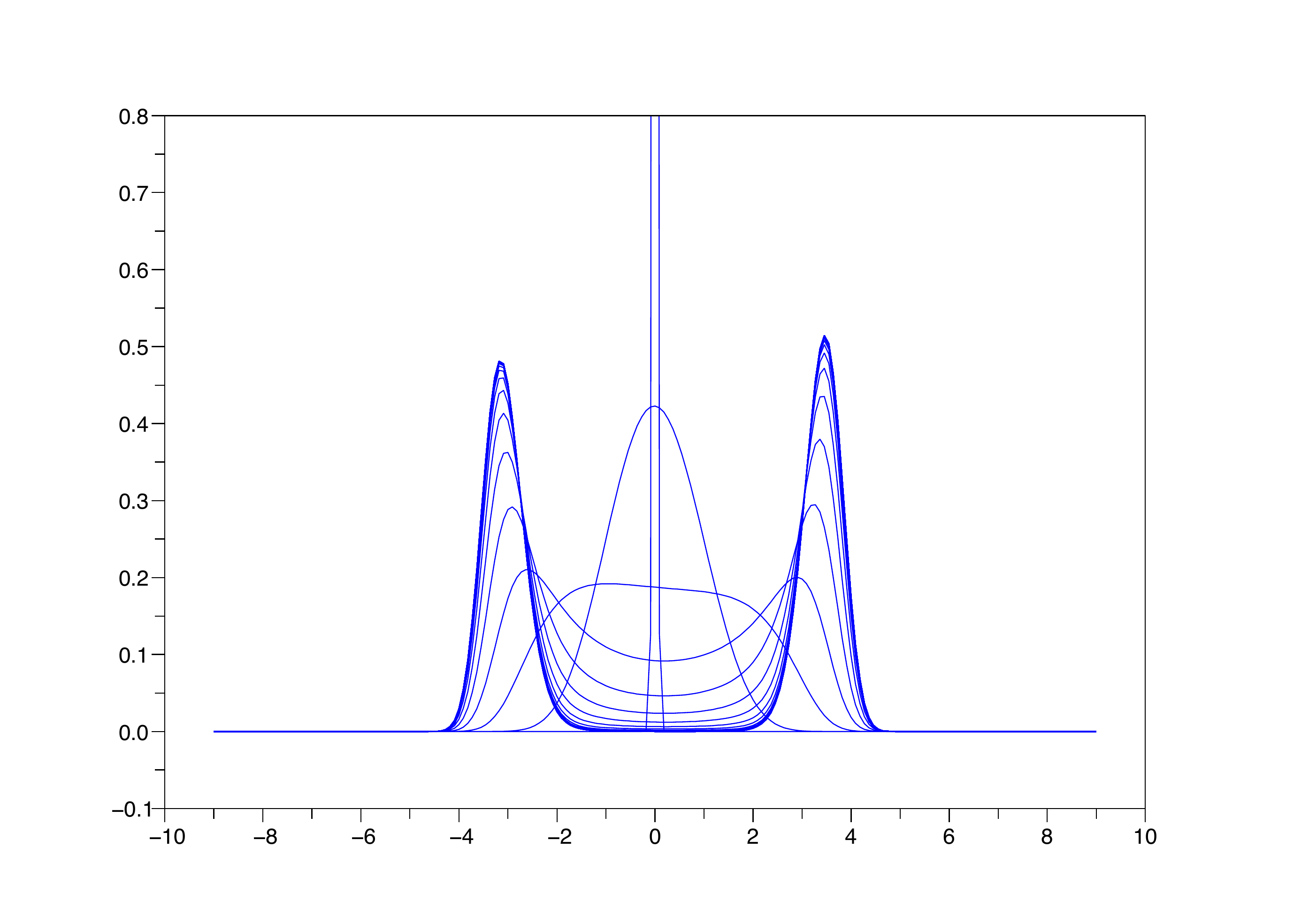}\
\includegraphics[width=6cm, angle=-0]{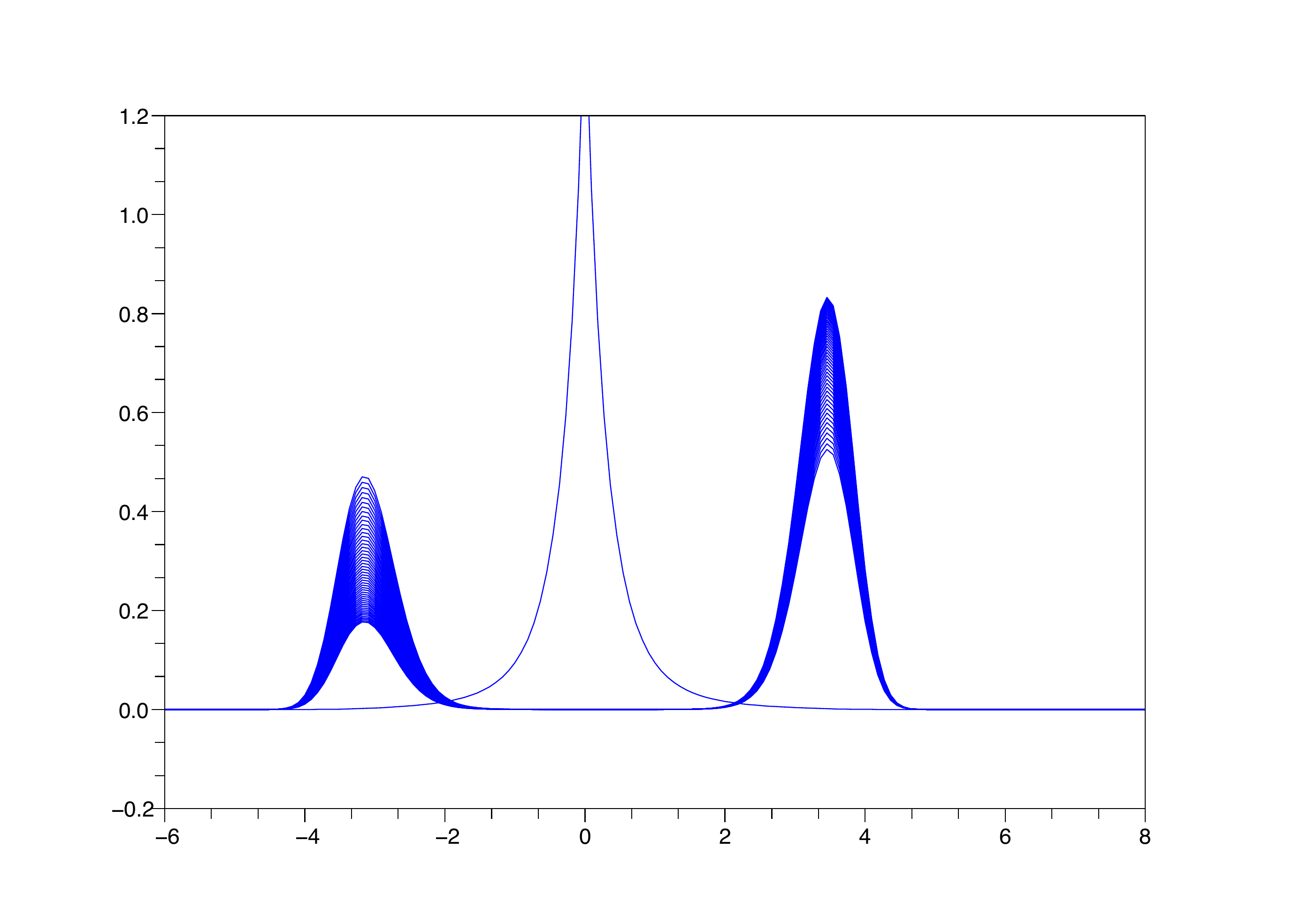}
} \caption{Evolution in time of the distribution in the $y$
variable, for the  biasing parameter $\Delta \lambda=  0.01$ with
snapshot every dt= 20 or 0.2 s (left) and  $20^4$ or 200s
(right).} \label{evolution_y}
\end{figure}

We are here interested also on the time behavior of the solution
$q(y,t)$ of the 1D Fokker-Planck equation \eqref{FPy}. For
instance, we may compare the time evolution of momentum for the 2D
and the reduced 1D problem. Thus, we need to compute not only the
stationary solution of equation \eqref{FPy}, but also its time
dependent solution. We choose to discretize equation \eqref{FPy}
using implicit in time finite difference method. The evolution of
the 1D reduced model illustrates again the slow-fast character of
this problem. In fact, we observe in Figure \ref{evolution_y} the
evolution in time of the density $q(y,t)$ for small (left) and for
large (right) times respectively. The convergence toward the final
stationary state is quite slow compared to the fast division
toward the two bump distribution at the initial stages.

We describe now how to recover all the moments of the partial
distribution function $p_\varepsilon$ in the $(\nu_1, \nu_2)$
plane, using the probability distribution function $q(y)$ solution
of \eqref{FPy} and the approximated slow manifold $x^*(y)$.

The function $p_\varepsilon$ is concentrated along the the curve
$\nu=(\nu_1(y), \nu_2(y))$ given by $\nu^{eq} + P(x^*(y),y)^T$,
see \eqref{chvar}. We parametrized this curve by means of a
curvilinear coordinate and define
$$
V(y) =  \| P(x'^*(y),1)^T \| \, .
$$
Then, for any test function $\Psi=\Psi(\nu_1,\nu_2)$, the moment
$M_\Psi$ of the probability distribution function
$p_\varepsilon=p_\varepsilon(\nu_1,\nu_2)$ is defined by
$$
M_\Psi = \int_\Omega \Psi \, p_\varepsilon\, d\nu_1 d\nu_2 \, ,
$$
and given by
$$
M_\Psi := \int \Psi( (\nu_1(y),\nu_2(y)) q(y) dy \,.
$$
This formulae can be used to compute either classical moments of
$p_\varepsilon$ or marginals by choosing e.g.
$\Psi=\delta_{\{\nu_1= \mu\}}$ to get the $\nu_2$-marginal as a
function of $\mu$. Note that $q$ has to be normalized in such way
that its total mass (along the slow manifold) is equal to $1$ i.e.
$M_{\Psi \equiv 1} =1$.

\begin{figure}[H]
\centering{\includegraphics[width=6cm,
angle=0]{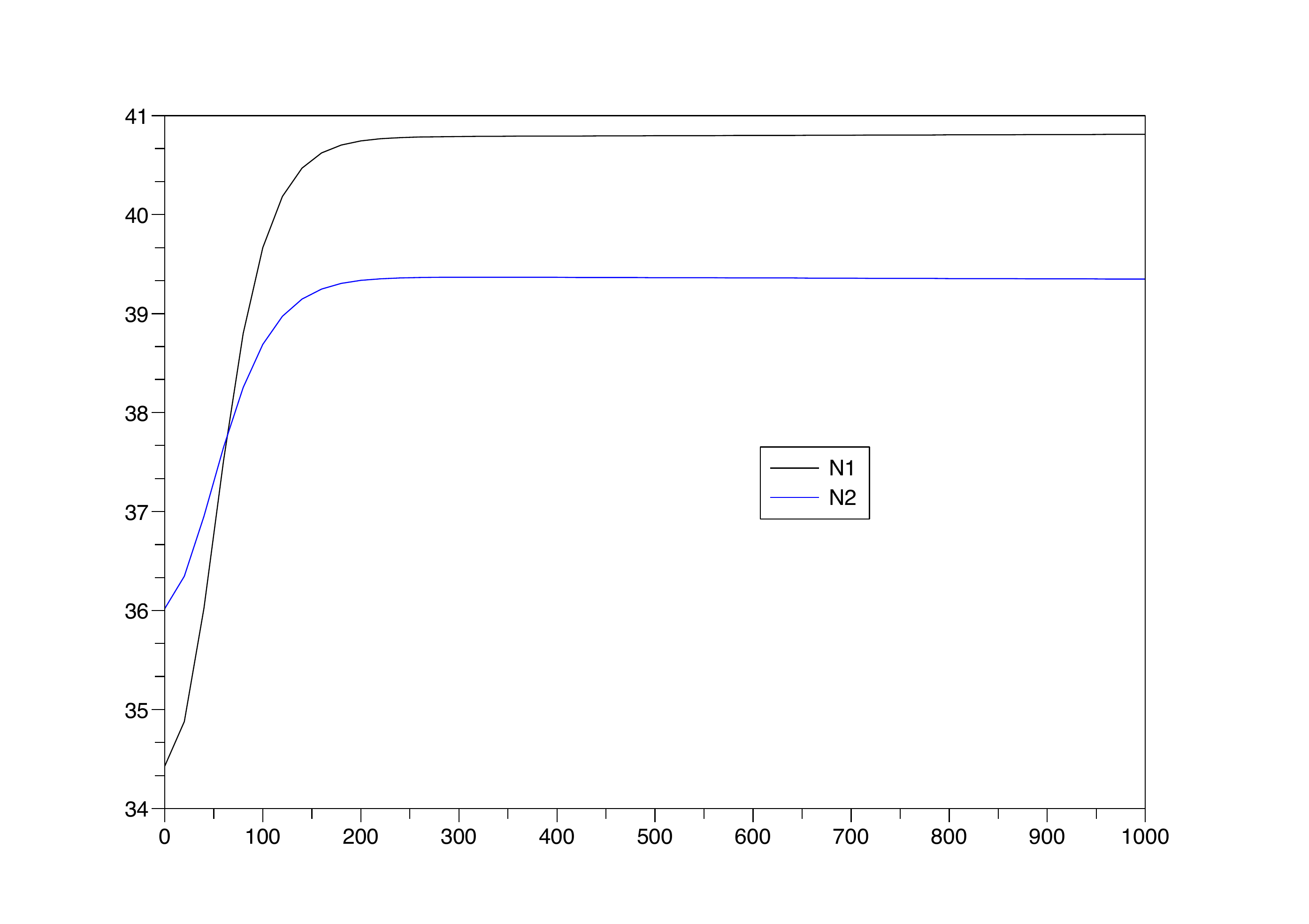} \
\includegraphics[width=6cm, angle=0]{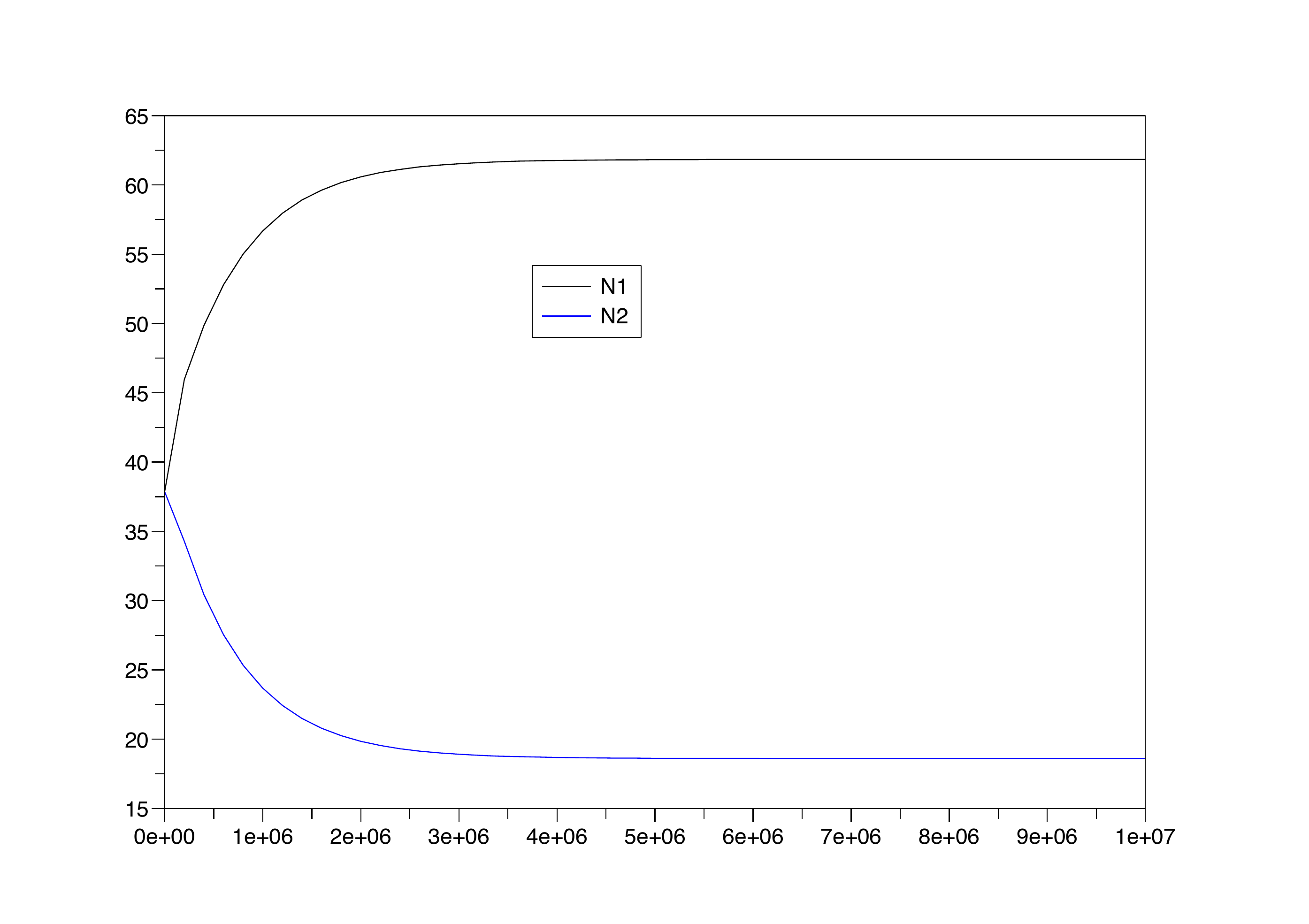}}
\caption{Evolution in time of the first order moments for the 1D
Fokker-Planck reduction. Left Figure: final time is $10^3$. Right
Figure: final time is $10^7$. Time unit is $0.01s$}
\label{1Dslowfast}
\end{figure}

Let us illustrate this metastability by the evolution of the first
moments of the distribution in Figure \ref{1Dslowfast}. The
initial data is a Dirac measure located above the spontaneous
point ($x=0$, $y>0$ small). We choose $\Delta \lambda=0$,
$\beta=0.3$ and $\omega_+=2.35$. We use a implicit scheme in order
to have no stability constraint on the time step. The number of
discretization point is 200 and the time step is $\Delta t=0.01$
for the left plot. It shows the fast dynamics: first the Dirac
measure diffuses onto a Gaussian blob and moves quickly toward the
spontaneous (unstable) state, then the Gaussian blob splits in the
two bumps around the two stable equilibrium points. It seems that
the solution has reached a equilibrium but it evolves very slowly.
The figure on the right corresponds to $\Delta t= 100$ and shows
this slow evolution toward the real equilibrium state. We will
comment more below.

\begin{figure}[H]
\centering{
\includegraphics[width=6cm, angle=0]{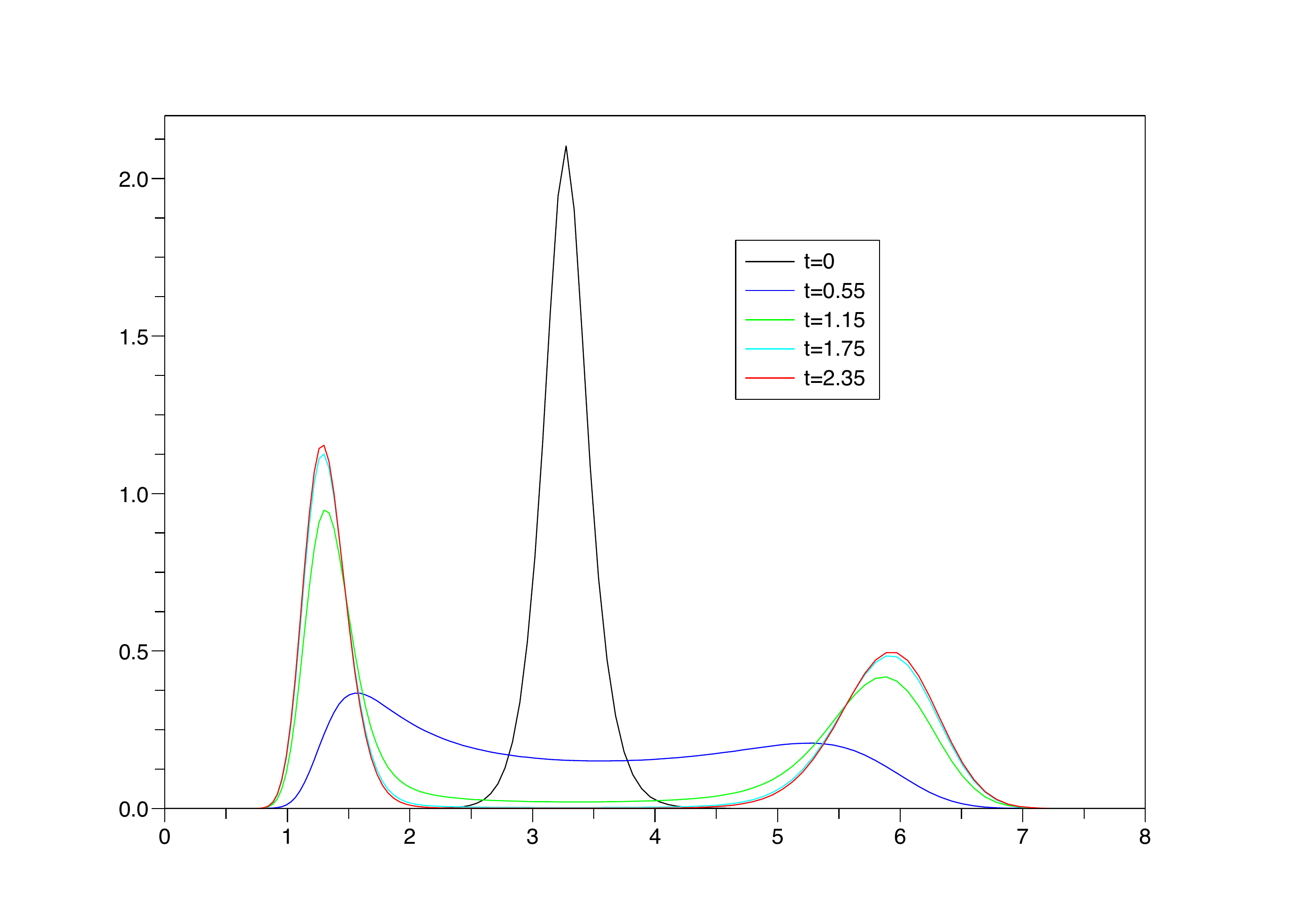}\
\includegraphics[width=6cm, angle=-0]{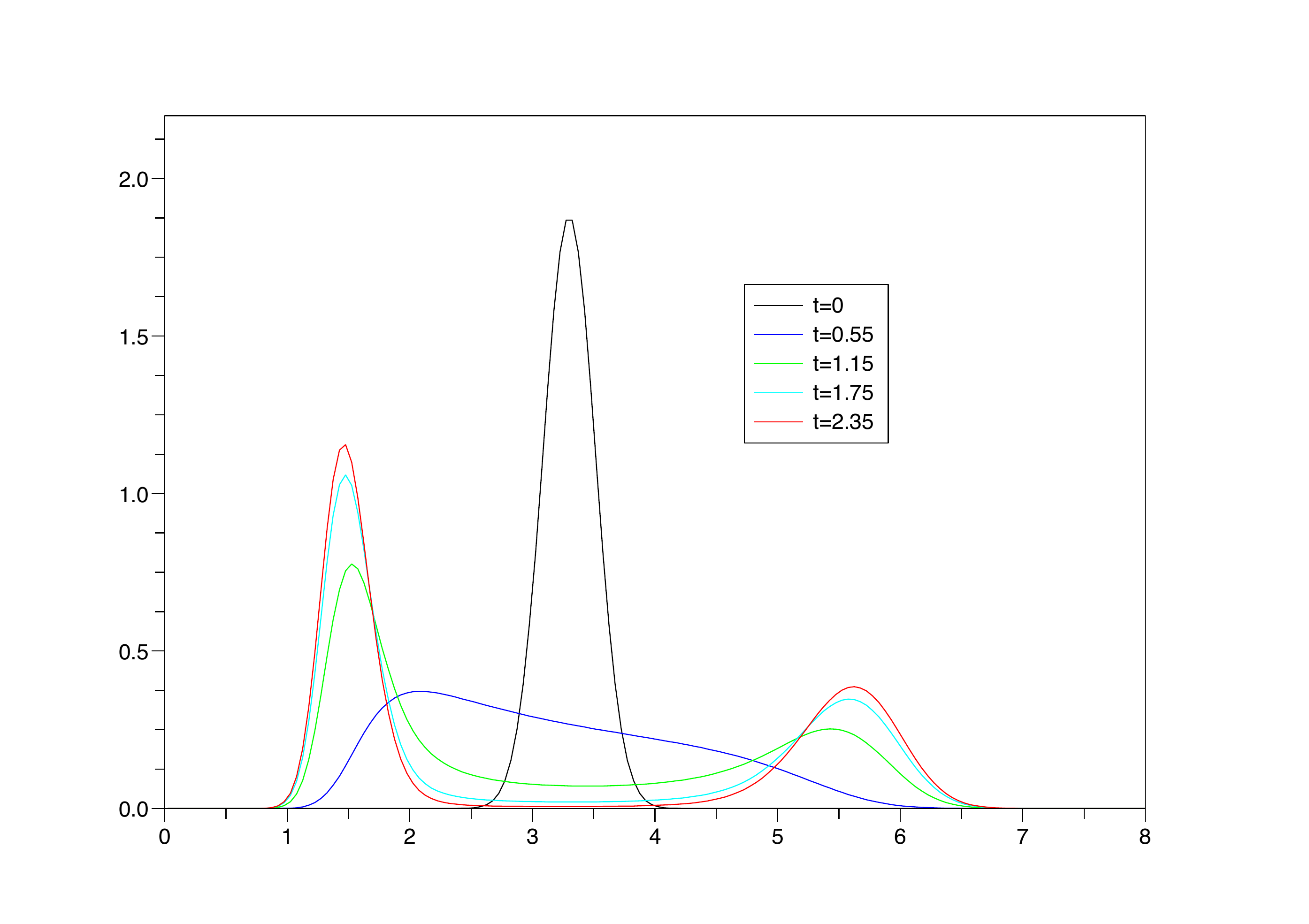}
} \caption{Comparison of the marginals in variable $\nu_1$, for
the  biasing parameter $\Delta \lambda=  0.01$. Left Figure: 1D
reduced model. Right Figure: 2D simulation.}
\label{comparisonmarg}
\end{figure}

We can finally compare the marginals in $\nu_1$ for the 1D reduced
model and the 2D simulations in Figure \ref{comparisonmarg}. We
can conclude that the transients of the 2D are captured extremely
well by the 1D reduced model.

\section{Reaction time and Performance}\label{probability}

In the previous sections, we have numerically studied the accuracy
of the reduced 1D model with respect to the 2D original one. We
discuss now some other information we can obtain from the 1D
problem, namely: the escaping time, section \ref{escaping}, and
the probability density to belong to a sub-domain of the phase
space, section \ref{section-rho}.

\subsection{Escaping time}\label{escaping}
\begin{figure}[H]
\centering{\includegraphics[width=10cm]{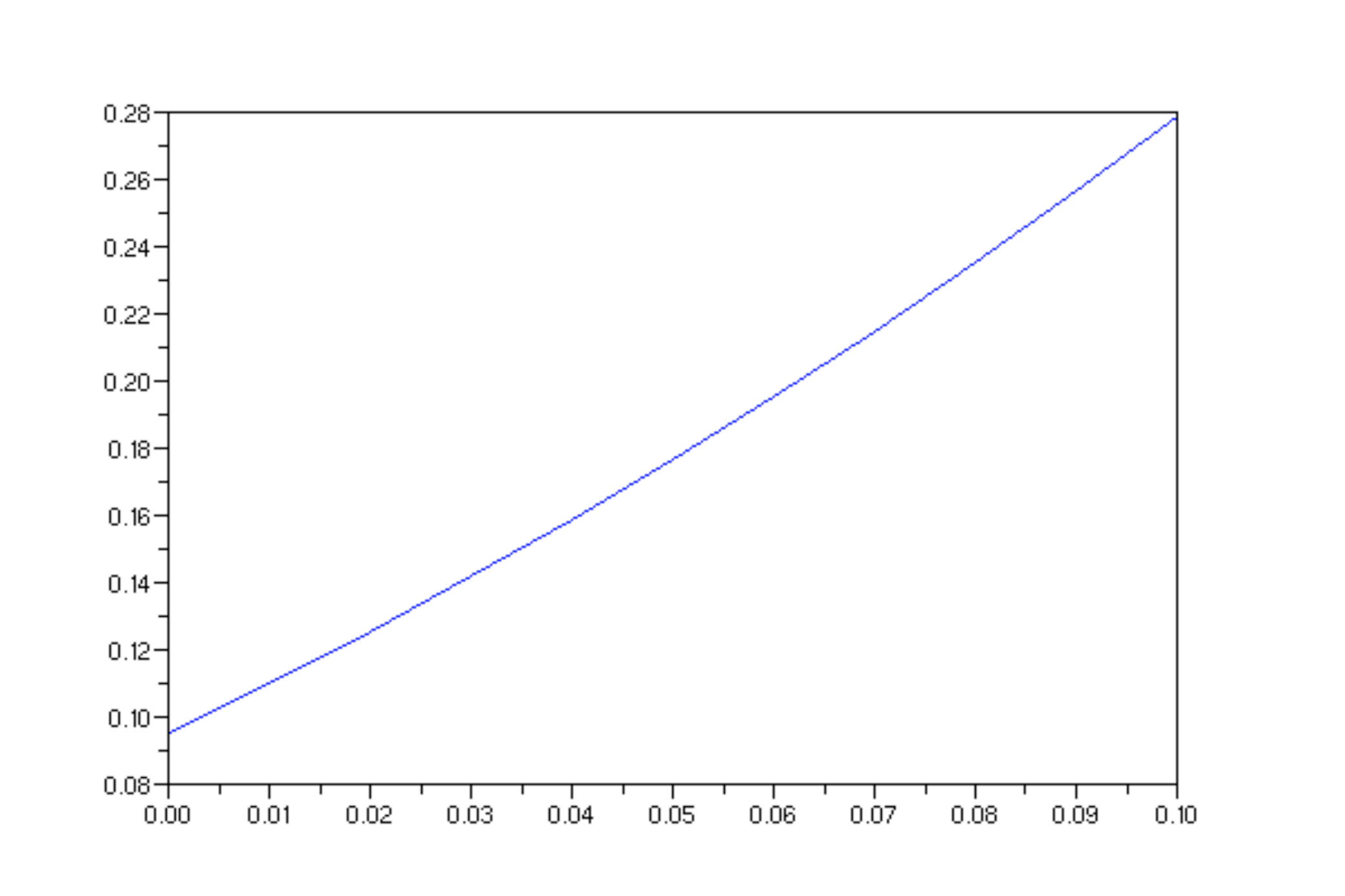}}
\caption{Computed potential gap as a function of $\Delta \lambda$.
} \label{DGG}
\end{figure}

Fixed a bias $\Delta \lambda$ and for a variable $\beta$ we can
easily compute the escaping time. In fact we recall that the
Kramers formula \cite{Gar}:
\begin{equation*}%\label{kramer}
\mathbb{E}(t)\sim \exp(H_G/\beta^2)
\end{equation*}
where $H_G$ is the maximum difference of the potential $G$
\begin{equation}\label{Hgap}
H_G = G_{max}-G_{min},
\end{equation}
apply in the one dimensional framework, without needing to compute
the solution $q(t,y)$ of the Fokker-Planck equation \eqref{FPy}.
Recall that $G_{max}$ corresponds to the potential value at the
spontaneous state while $G_{min}$ corresponds to the minimum of
the potential value at the two decision states. In Figure
\ref{DGG}, we plot the potential gap computed by means of
\eqref{Hgap} as a function of the bias $\Delta \lambda$.

In the 2D problem, since the drift is not the gradient of a
potential, Kramer's rule does not apply and the escaping times can
only be computed for the unbiased case, $\Delta \lambda=0$. In
fact, for $\Delta \lambda=0$ the problem is symmetric in $\nu_1$
and $\nu_2$ and thus, we know that the firing rates will separate
in two identical bumps. Then, starting the computation from an
initial data narrowly concentrated around one stable equilibrium
point (say $S_1$), the escaping time $T$ can be defined as the
time needed to have half of the total mass moving to a
neighborhood of $S_3$. In particular, the expectation
$\mathbb{E}(T)$ has an exponential behavior and its associated
potential gap is $H_G = 0.1$. Of course, this kind of argument
cannot be extended to the biased case and thus, the 1D reduced
model is essential.

\subsection{Probability densities - Performance }\label{section-rho}

We can compare the value of the probability $\rho_i$ of the firing
rates to be in some domain $\Omega_i$ for the 2D Fokker-Plank
model and the 1D reduced FP model. In particular, we shall
compare: $\rho_p$ the probability for the 2D problem that at time
$t$ the firing rates belong to $\{(x,y): y>0\}$; with $\rho_+$
representing the probability for the solution of the 1D problem to
belong to $\Omega_+=[0,y_m]$, see Figure \ref{chang-variable}. We
fix the standard derivation $\beta=0.1$, and let the bias $\Delta
\lambda$ varying from $0$ to $0.05$, since the values of $\rho_p$
and $\rho_+$ for $\Delta \lambda$ bigger than $0.03$ are already
very close (the relative error being of the order of $10^{-4}$).

\begin{figure}[H]
\centering{\includegraphics[width=7cm,
angle=-90]{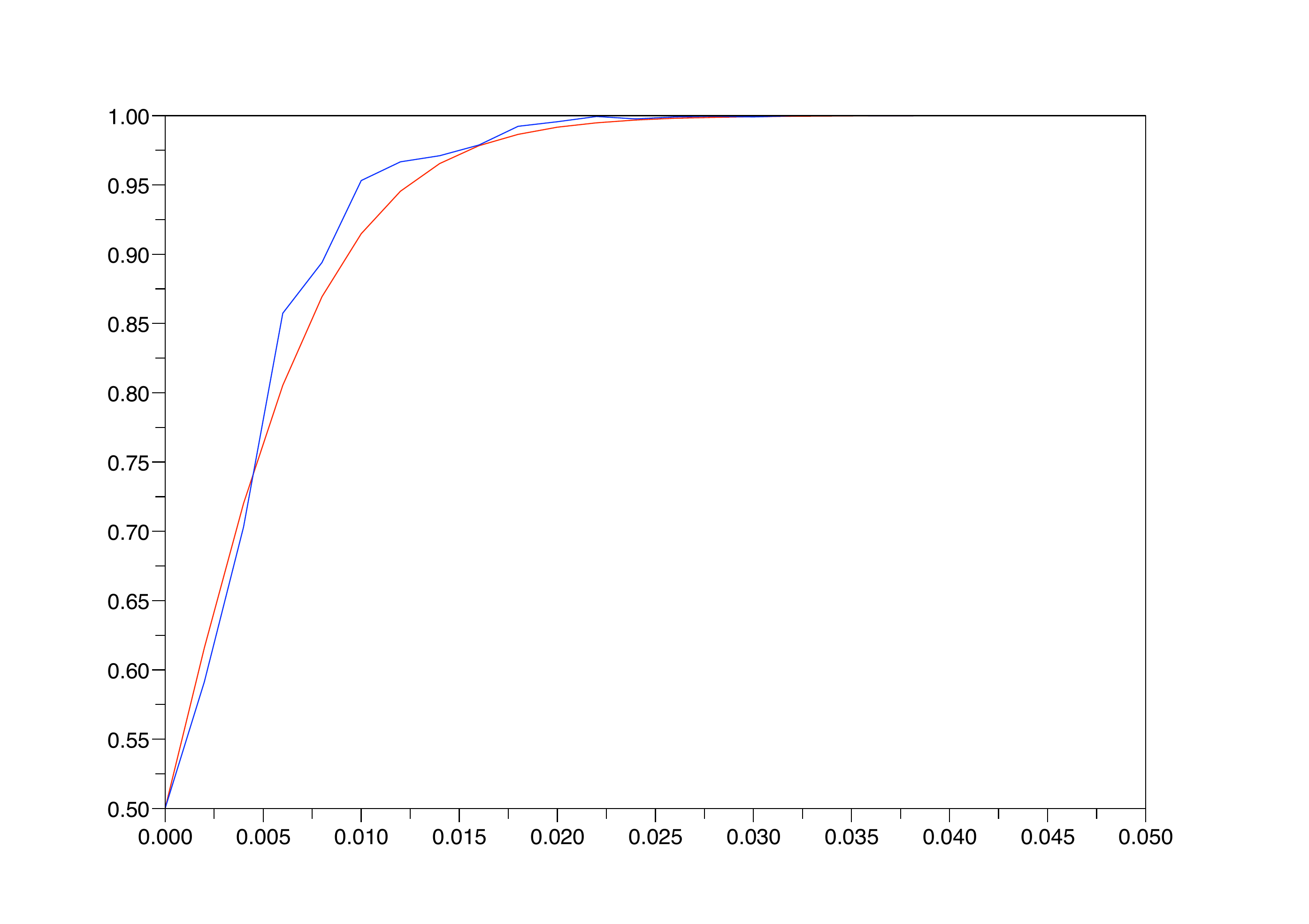}} \caption{Comparison of the values
for $\rho_+$ and $\rho_\infty$ with respect  to the biasing
parameter $\Delta \lambda\in [0,0.05]$. Red line: the values
computed from the 1D reduced problem. Blue line: the values
extrapolated from the 2D problem.} \label{rhop}
\end{figure}

We recall that, in the 2D  problem, we has to wait for a very long
time in order to reach equilibrium, since we have a meta-stable
situation, see \cite{CCM}. Nevertheless, we note that the
$\rho_p(t)$ profile is exponentially increasing converging to an
asymptotic value $\rho_\infty$. We can extrapolate this value from
the values of $\rho_p(t)$ for some initial iterations as follows.
Assume that the probability $\rho_p$ behaves like:
\begin{equation*}%\label{rho1}
\rho(t)=\rho_\infty - a \exp(-t/\tau)
\end{equation*}
where $\rho_\infty$, $a$ and $\tau$ have to be determined  by  an
``exponential regression''. For a sequence of time $t_i$ (of the
form $t_i=t_0+ T*i$, $i=0\cdots N$ that corresponds to the
computed values of $\rho_p$), we define $\Delta \rho_i$ as the
difference $\rho_p(t_i)-\rho_p(t_i+T)$, we get:
\begin{equation}\label{delta_rho}
\Delta \rho_i = a \exp(-t_i/\tau)(1-\exp(-T /\tau)).
\end{equation}
Taking the $\log$ and the difference between two indexes $i$ and $j$
we obtain the expected value of $\tau$ as:
\begin{equation*}%\label{tau}
\tau \approx - {\log(\Delta \rho_i) - \log(\Delta \rho_j) \over  t_i - t_j}.
\end{equation*}
Finally from \eqref{delta_rho}, knowing $\tau$ (and $T$), we
recover $a$, and the asymptotic limit $\rho_\infty$ is uniquely
determined by:
\begin{equation}\label{rho-infinity}
\rho_\infty = \rho(t_0) + a \exp(-t_0/\tau).
\end{equation}
We show in Figure \ref{rhop} the comparison between the values for
the one dimensional computation (red line) and the one
extrapolated from the 2D computation, see equation
\eqref{rho-infinity}, using a final time of 20 seconds (blue
line). Note that the non-smoothness of the blue line (2D
extrapolation) may be due to the fact that for computing $\rho_p$
we need to compute the inner product for any point $\nu$ of the
phase space $(\nu_1, \nu_2)$ : $<(\nu_{eq} -\nu), P_{1,j}> $, and
we choose for different values of $\Delta \lambda$ the same
equilibrium point $\nu_{eq}$ and matrix $P$: for instance, for
$\Delta \lambda= 0.016, 0.018, 0.020, 0.022, 0.024$ we choose the
values of $\Delta \lambda = 0.020$:
$$\nu_{eq}=\left(
\begin{array}{c}
3.0448158\\
3.2397474
\end{array}
\right) \; , \qquad
P=\left(
\begin{array}{cc}
0.7003255 & -0.6959201 \\
0.7138236  & 0.7181192
\end{array}
\right).$$

%%%%%%%%%%%%%%%%%%%%%%%%%%%%%%%%%%%%

\subsection*{Acknowledgments}
J.~A.~Carrillo is partially supported by the projects
MTM2011-27739-C04 from DGI-MICINN (Spain) and 2009-SGR-345 from
AGAUR-Generalitat de Catalunya. S.~Cordier and S.~Mancini are
partially supported by the ANR projects: MANDy, Mathematical
Analysis of Neuronal Dynamics, ANR-09-BLAN-0008-01. All authors
acknowledge partial support of CBDif-Fr, Collective behaviour $\&$
diffusion : mathematical models and simulations,
ANR-08-BLAN-0333-01.

%%%%%%%%%%%%%%%%%%%%%%%%%%%%%%%%%%%%

\end{document}